# LARGE AND MODERATE DEVIATIONS IN POISSON NAVIGATIONS

PARTHA PRATIM GHOSH ,* ** *Technische Universität Braunschweig*
BENEDIKT JAHNEL,* *** *Technische Universität Braunschweig & Weierstrass Institute Berlin*
SANJOY KUMAR JHAWAR,**** *INRIA Paris & Telecom Paris*

## Abstract

We derive large- and moderate-deviation results in random networks given as planar directed navigations on homogeneous Poisson point processes. In this non-Markovian routing scheme, starting from the origin, at each consecutive step a Poisson point is joined by an edge to its nearest Poisson point to the right within a cone. We establish precise exponential rates of decay for the probability that the vertical displacement of the random path is unexpectedly large. The proofs rest on controlling the dependencies of the individual steps and the randomness in the horizontal displacement as well as renewal-process arguments.

## 1. Introduction

We consider decentralized traffic-flow networks or drainage networks, where nodes transmit data individually to dedicated neighboring nodes according to local rules, or respectively liquid flows towards a certain direction. For this, directed navigations on homogeneous Poisson point processes in the Euclidean plane can be used as an underlying network model [3, 5]. The traffic generated at each node is forwarded iteratively to the nearest node, for example, to the right in the horizontal direction. While analyzing such models, it has been observed in [10], that the spatial traffic-flow density at a given location in the domain, i.e. the spatial average of accumulated traffic at a microscopic area around that point, follows a law of large numbers in the high-density limit. Here, the asymptotic traffic density can be captured as an integral in terms of the spatial intensity of the point process along with their rate of traffic generation. The bulk contribution in the limit in this asymptotics comes from the traffic generated within a

Received 1 July 2025; accepted 1 July 2025.

* Postal address: Technische Universität Braunschweig, Universitätsplatz 2, 38106 Braunschweig, Germany.
** Email address: p.pratim.10.93@gmail.com
*** Postal address: Weierstrass Institute Berlin, Mohrenstrasse 39, 10117 Berlin, Germany. Email: benedikt.jahnel@tu-braunschweig.de
**** Postal address: INRIA Paris, 48 Rue Barrault, 75013 Paris, France. Email: sanjoy-kumar.jhawar@inria.fr

thin horizontal strip at that given point due to sub-ballisticity of the trajectory of the navigation. This sub-ballistic behavior is nothing but a version of the straightness condition studied in [11].

Going beyond sub-ballisticity and the law of large numbers, questions about the deviation behavior of paths in the traffic network arise. In a stationary setting, is it rare to have a path originating from the origin to deviate much away from the horizontal axis? For this, we denoted by $\mathcal{Y}_t$ the vertical position of the path measured as a function of horizontal position $t$. If the answer to the question is no then essentially these paths are not so unlikely to miss the microscopic region around a given location, and the traffic accumulated along the path may not contribute to the aggregated traffic at the location of interest. This in turn would lead to deviations in the throughput at the target location and it is an indication of unexpected network performance. On the other hand, if the answer is yes, that means that the paths do not deviate too much and contribute to the bulk of the throughput. But how rare is the deviation event? In other words, at what scale does the probability of such an event live on and what are the rates? To answer these sort of questions, in this paper we investigate the large- and moderate-deviation behavior of the deviation of these paths away from the horizontal axis in a Poisson navigation model where the consecutive *successors* are defined as the closest neighbor within a cone to the right of the current point in the Euclidean plane.

Large- and moderate-deviation principles constitute a cornerstone in the analysis of random geometric models and, in general, in stochastic processes in the last four to five decades, well documented, for example, in [7, 8, 17]. These principles essentially capture rare events by quantifying their asymptotic unlikeliness. The small probability of the events is expressed in terms of an exponential rate and a rate function. In the context of telecommunication networks, large-deviation principles can be used to point out how badly the network is performing and what the most prominent reasons for this behavior are. The large-deviation principle has been studied for many different telecommunication network models, mainly identifying a poor service indicator; see, for example, [12, Chapter 6] and references therein. On the other hand, moderate deviations are nothing but large deviations with a slower scaling. The moderate-deviation principle for a sequence of independent and identically distributed (i.i.d.) centered random variables has been studied, for example, in [9] as a much cleaner and shorter version of [1] and for sequences of Banach-space-valued random variables, see the classic work of Ledoux [13]. In the case of moderate deviations, the rate function turns out to be of Gaussian type irrespective of the model under consideration, whereas in large deviations the rate function depends on characteristics of the random variables.

The model we study in this paper is the underlying setup for the traffic-flow network mentioned above. It is known as a *navigation* on Poisson point processes in Euclidean space and specifically the directed and radial navigation were introduced in [5]. However, we introduce an additional parameter $0 < \theta \leq \pi/2$ that controls the angle of the cone in which the navigation searches for the next neighbor. By definition, if the navigation starts at multiple points, this type of navigation gives rise to a tree structure. In applications to telecommunication networks, the root of the spanning trees can be seen as a network head where all the information is gathered along the edges of the tree and processed. Hence, it is worthwhile to know the local and global structural properties of the tree. In the work of Baccelli and Bordenave [3], which only deals with the case $\theta = \pi/2$, the quantities of interest were the local tree functionals around a vertex, for example, the degree, properties of the path in the tree from a vertex to the root, namely the total length or properties of the tree structure, and its shape. In a somewhat more geometric and analytical work, the convergence of the tree to a Brownian web under the appropriate scaling limit has been studied in [6, 15] for the Poisson setting, in general dimensions, and more recently, for a discrete setting on perturbed planar lattices [16]; however, all this is only in the case $\theta = \pi/2$. In this paper, in particular, we improve the polynomial-decay properties in [3, Lemma 4.11] to stretched-exponential decay of a moderate-deviation type as well as go beyond the associated central limit theorem [3, Theorem 4.7]. Let us also mention the work of Bonichon and Marckert [4], where navigation characteristics are investigated in a high-density regime for not necessarily homogeneous Poisson point processes.

When $\pi/4 < \theta \leq \pi/2$, the model of navigation possesses a challenging dependence structure along its path, since each step of the navigation carries an extra piece of information from the previous steps via nontrivial regions that have been observed as void spaces for the underlying Poisson point process. We call this piece of information the *history set* and identify the steps where the history set is well behaved, giving rise to a renewal structure that is essential for our analysis. However, the tail properties of the step variables of the renewal process are challenging to analyze, which makes up for the bulk of the technical part of the paper. For this, we control the exponential decay of the inter-stopping time gaps via bounds on the dynamics of the width of the history set, using Markov-chain comparison ideas from [6, Proposition 3.1]. However, this control fails to be sufficiently detailed in order to establish the large deviations in the whole parameter regime. Indeed, the model becomes substantially simpler for $\theta \in (0, \pi/4]$ as the steps become i.i.d. and we can then provide the full large- and moderate-deviation analysis with rather explicit rate functions.

Let us mention that our main results should still hold in higher dimensions with appropriately adjusted definitions, for example, of higher-dimensional angles. The associated proofs should not pose any additional substantial difficulties. In order to keep the exposition more accessible, we focus here only on the planar case. Similarly, the qualitative picture should remain intact if the directed navigation defined below, which is based on a nearest-neighbor relation, is replaced by a navigation based on connecting, for example, to the $k$th nearest neighbor, $k \geq 1$. However, at least for $\theta \in (\pi/4, \pi/2]$, many of the underlying definitions and arguments would need to be changed, but a corresponding renewal structure should still be present. Another interesting direction for future research could be to assign individual angles to every vertex via i.i.d. marks and study the resulting navigation process.

Organization of the paper is as follows. In Section 2 we describe the model in detail, state the main results about large- and moderate-deviation principles and the scaling property of the rate functions with respect to the intensity of the Poisson point process. In Section 3 we first uncover the hidden independence structure in the model and state the key supporting results that enable us to prove the main result for the moderate deviations for $\theta < \pi/2$. In Section 4 we discuss the large-deviation principle in the dependent case, i.e. for $\pi/4 < \theta < \pi/2$, and prove it under a key assumption for the tail behavior of the renewal-step variables. Section 5 contains the proofs of all the supporting results that are stated in Section 3. Section 6 contains the separate argument for the moderate-deviation principle in case $\theta = \pi/2$ and in Section 7 we prove all statements regarding large deviations. Finally, in the Appendix A we elaborate on a connection between moderate and large deviations.

## 2. Setting and Main Results

Let $\mathcal{P}_\lambda$ denote the homogeneous Poisson point process with intensity $\lambda > 0$ on $\mathbb{R}^2$ and consider an additional point $o$ at the origin such that $\mathcal{P}_\lambda \cup \{o\}$ is the Poisson point process under its Palm distribution $\mathbb{P}$. We are interested in the large and moderate deviations of a dependent sequence of waypoints, starting at $o$, that are defined iteratively by choosing a successor point, which is the closest Poisson point towards the right within a cone. More precisely, for any $v \in \mathbb{R}^2$, let $(r_v, \varphi_v) \in \mathbb{R}_+ \times [-\pi, \pi)$ denote its polar coordinates, with the first unit vector $e_1$ corresponding to the polar coordinates (1,0), and consider

$$\mathcal{C}_\theta := \{v = (r_v, \varphi_v) \colon r_v > 0, |\varphi_v| \leq \theta\},$$

the cone with angle $2\theta$ centered at the origin. For $0 \leq \theta \leq \pi/2$, we might say that the cone is facing towards the right. For any $v \in \mathbb{R}^2$, we write $\mathcal{C}_\theta(v) := v + \mathcal{C}_\theta$ and $\mathcal{C}^o_\theta(v)$ for the interior of $\mathcal{C}_\theta(v)$. Then, we consider the following family of navigations based on the usual Euclidean metric $|\cdot|$ in $\mathbb{R}^2$.

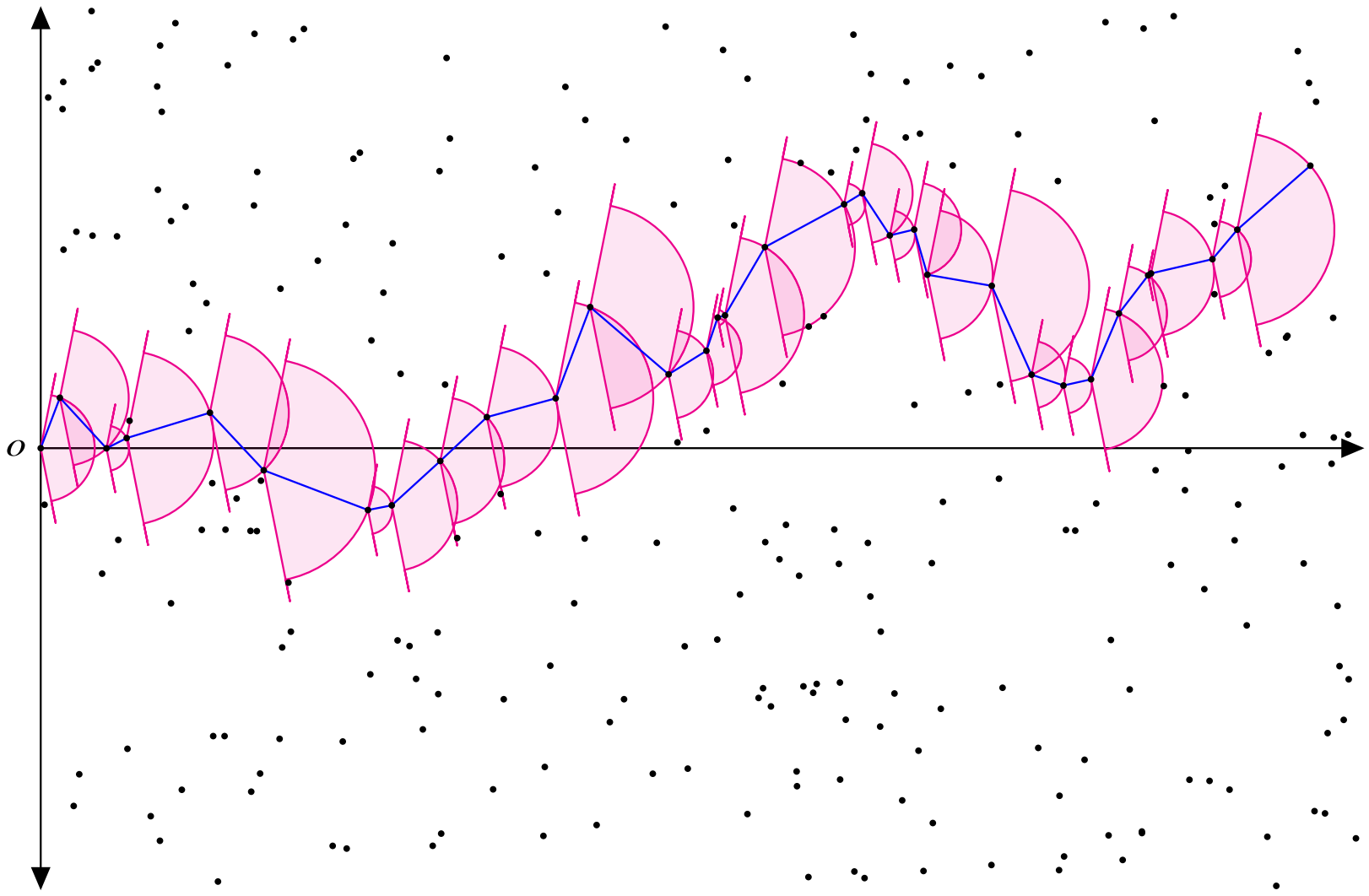

FIGURE 1. Simulated sample path of $\bar{\mathcal{V}}$ for $\theta = \arctan(5)$.

**Definition 1.** (*Navigations.*) Let $V_0 = o$. Then, we call $\mathcal{V} := \{V_n\}_{n\geq 0} \subset \mathcal{P}_\lambda$, iteratively defined as

$$V_{i+1} := \operatorname{argmin}\{|v - V_i| \colon v \in \mathcal{P}_\lambda \cap \mathcal{C}_\theta(V_i)\},$$

a directed $\theta$ navigation. Here $V_{i+1}$ is denoted the successor of $V_i \in \mathcal{P}_\lambda$ and $U_i := V_i - V_{i-1}$ the ith progress of the navigation.

Note that, under $\mathbb{P}$, the argmin is uniquely defined almost surely (a.s.). We are interested in a continuous-time process based on the navigation $\mathcal{V}$. For this, consider $\bar{\mathcal{V}} := \bigcup_{k\geq 0} [V_k, V_{k+1}]$, the *interpolated trajectory* for the navigation, where $[V_k, V_{k+1}] \subset \mathbb{R}^2$ should be understood as the one-dimensional line segment that connects $V_k$ and $V_{k+1}$. Then, we can see the interpolated trajectory $\bar{\mathcal{V}}$ as a piecewise affine and continuous path, parametrized with respect to time $t$ as $\{(t, \mathcal{Y}_t)\}_{t\geq 0}$, where the parameter $t$ denotes the progress along the $x$ axis and $\mathcal{Y}_t \in \mathbb{R}$ is the corresponding $y$ coordinate or *vertical displacement* at time $t$. In particular, for $t = \pi_1(V_k)$, where $\pi_1$ denotes the projection to the first Cartesian coordinate, we have $\mathcal{Y}_t = \pi_2(V_k)$, for every $k \geq 0$, where $\pi_2$ denotes the projection to the second Cartesian coordinate; see Figure 1 for an illustration.

The following standard notation is used throughout the paper. For any set $\Gamma$, we denote its interior by $\Gamma^o$, its closure by $\overline{\Gamma}$, and its complement by $\Gamma^c$. Our first main result is the *moderate-deviation principle* for $\mathcal{Y} = \{\mathcal{Y}_t\}_{t\geq 0}$.

**Theorem 1.** *For any* $0 < \lambda$, $0 < \theta \leq \pi/2$ *and* $0 < \varepsilon < 1/2$, $\{t^{-1/2-\varepsilon}\mathcal{Y}_t\}_{t\geq 0}$ *obeys the moderate-deviation principle with rate* $t^{2\varepsilon}$ *and rate function* $I_{\lambda,\theta}(x) := \rho(\lambda, \theta)x^2$, *where* $\rho(\lambda, \theta) > 0$, *meaning that, for any Borel set* $\Gamma \subseteq \mathbb{R}$,

$$
\begin{aligned}
-\inf_{x\in\Gamma^o} I_{\lambda,\theta}(x) &\leq \liminf_{t\to\infty} t^{-2\varepsilon} \log \mathbb{P}(t^{-1/2-\varepsilon}\mathcal{Y}_t \in \Gamma) \\
&\leq \limsup_{t\to\infty} t^{-2\varepsilon} \log \mathbb{P}(t^{-1/2-\varepsilon}\mathcal{Y}_t \in \Gamma) \\
&\leq -\inf_{x\in\overline{\Gamma}} I_{\lambda,\theta}(x).
\end{aligned}
$$

*Additionally,* $\rho(\lambda, \theta)$ *satisfies the scaling relation* $\rho(\lambda, \theta) = \sqrt{\lambda}\rho(1, \theta)$.

The constant $\rho(\lambda, \theta)$ can be generally expressed in a semi-explicit way; see (6) in Section 3. In the case $0 < \theta \leq \pi/4$, it is given by

$$
\rho(\lambda, \theta) := \frac{\int_0^\infty \mathrm{d}r\, r^2 \exp(-\lambda\theta r^2) \int_{-\theta}^{\theta} \mathrm{d}\varphi\, \cos\varphi}{2\int_0^\infty \mathrm{d}r\, r^3 \exp(-\lambda\theta r^2) \int_{-\theta}^{\theta} \mathrm{d}\varphi\, \sin^2\varphi} = \frac{\sqrt{\pi\lambda\theta}\,\sin\theta}{2\theta - \sin(2\theta)}. \tag{1}
$$

Let us mention that Theorem 1 can readily be used to establish also the strong law of large numbers $t^{-1}\mathcal{Y}_t \to 0$ a.s. In fact, Theorem 1 implies that even $t^{-1/2-\varepsilon}\mathcal{Y}_t \to 0$ a.s. for all $0 < \varepsilon < 1/2$.

For all sufficiently small angles, we show that $\mathcal{Y}$ also satisfies the *large-deviation principle*. The rate function is given in terms of the multivariate rate function of i.i.d. progress variables defined in terms of the usual scalar product $\langle\cdot, \cdot\rangle$ in $\mathbb{R}^2$.

**Lemma 1.** *Consider i.i.d. copies* $\{\tilde{U}_i\}_{i\geq 1}$ *of the progress variable* $U_1 \in \mathbb{R}^2$ *in the directed* $\theta$ *navigation. Then,* $\{n^{-1}\sum_{i=1}^n \tilde{U}_i\}_{n\geq 1}$ *satisfies the large-deviation principle with rate n and rate function*

$$
\mathcal{J}_{\lambda,\theta}(u) := \sup\{\langle\gamma, u\rangle - J_{\lambda,\theta}(\gamma) \colon \gamma \in \mathbb{R}^2\},
$$

*where, for all* $\gamma \in \mathbb{R}^2$,

$$
J_{\lambda,\theta}(\gamma) := \log\left(\lambda \int_0^\infty dr\, r \exp(-\lambda\theta r^2) \int_{-\theta}^{\theta} d\varphi\, \exp(\gamma_1 r\cos\varphi + \gamma_2 r\sin\varphi)\right) < \infty.
$$

*Furthermore,* $u \mapsto \mathcal{J}_{\lambda,\theta}(u)$ *is continuous on* $\{u \in \mathbb{R}^2 \colon \mathcal{J}_{\lambda,\theta}(u) < \infty\} = \mathcal{C}_\theta^o$.

We are now in the position to state our second main theorem.

**Theorem 2.** *For any* $0 < \lambda$ *and* $0 < \theta \leq \pi/4$, $\{t^{-1}\mathcal{Y}_t\}_{t\geq 0}$ *obeys the large-deviation principle with rate t and rate function* $\mathcal{I}_{\lambda,\theta}(x) := \inf\{\beta\mathcal{J}_{\lambda,\theta}(1/\beta, x/\beta) \colon \beta > 0\}$, *meaning that, for any Borel set* $\Gamma \subseteq \mathbb{R}$,

$$
\begin{aligned}
-\inf_{x\in\Gamma^o}\mathcal{I}_{\lambda,\theta}(x) &\le \liminf_{t\to\infty} t^{-1}\log\mathbb{P}(t^{-1}\mathcal{Y}_t\in\Gamma)\\
&\le \limsup_{t\to\infty} t^{-1}\log\mathbb{P}(t^{-1}\mathcal{Y}_t\in\Gamma)\\
&\le -\inf_{x\in\overline{\Gamma}}\mathcal{I}_{\lambda,\theta}(x).
\end{aligned}
$$

*Additionally, $\mathcal{I}_{\lambda,\theta}$ satisfies the scaling relation $\mathcal{I}_{\lambda,\theta}=\sqrt{\lambda}\mathcal{I}_{1,\theta}$.*

We believe that $\mathcal{Y}$ also satisfies the *large-deviation principle* for $\pi/4<\theta<\pi/2$, but with a substantially more involved rate function $\mathcal{I}'_{\lambda,\theta}$ (see Theorem 3). As we highlight in the following sections, the challenge comes from the fact that, for $\pi/4<\theta<\pi/2$, the individual steps are not i.i.d. any more. The underlying renewal structure, mentioned in the introduction and explained in detail in the next section, which is also crucial for the proof of the moderate-deviation principle, has individual renewal steps with exponential tails. One renewal step therefore makes a negligible contribution on the moderate-deviation scale, but is non-negligible on the large-deviation scale. Due to the challenging dependence structure within the renewal steps, we are only able to give the large-deviation principle in Theorem 3 below assuming control on the tail behavior of the interarrival steps, see Assumption 1, and postpone the details to Section 4.

Before we enter a more detailed discussion about the proofs, let us elaborate on the connection between Theorems 1 and 2.

**Remark 1.** There is a heuristic that suggests that $\rho(\lambda,\theta)$ coincides with the second derivative of the large-deviation rate function at zero, that is, $\rho(\lambda,\theta)=\ddot{\mathcal{I}}_{\lambda,\theta}(0)/2$, which stems from the following intuition. With $h=xt^{\varepsilon-1/2}$ small, we roughly have, for large $t$,

$$
\begin{aligned}
t^{-2\varepsilon}\log\mathbb{P}(t^{-1/2-\varepsilon}\mathcal{Y}_t\approx x) &= t^{-2\varepsilon}\log\mathbb{P}(t^{-1}\mathcal{Y}_t\approx xt^{\varepsilon-1/2})\\
&\approx -t^{1-2\varepsilon}\mathcal{I}_{\lambda,\theta}(xt^{\varepsilon-1/2})\\
&= -x^2h^{-2}\mathcal{I}_{\lambda,\theta}(h)\\
&\approx -x^2\ddot{\mathcal{I}}_{\lambda,\theta}(0)/2,
\end{aligned}
$$

where we used the fact that $\mathcal{I}_{\lambda,\theta}(0)=\dot{\mathcal{I}}_{\lambda,\theta}(0)=0$ by the strong law of large numbers mentioned above. However, as we will investigate in detail in Appendix A, this intuition fails already in the independent case where $0<\theta\le\pi/4$. Roughly, this comes from the fact that, on the large-deviation level, there is a non-negligible interplay between the progress in the vertical and horizontal direction, which manifests itself also in $\ddot{\mathcal{I}}_{\lambda,\theta}(0)$. However, this is not the case on the level of moderate deviations. Here, deviations in the horizontal direction come at an exponential cost of rate $t$, and hence, play no role at rate $t^{2\varepsilon}$, with $0<\varepsilon<1/2$.

## 3. Strategy of Proof

The proof of Theorem 2 is a consequence of Lemma 1 and stepwise independence in the case $\theta\le\pi/4$. We present the arguments in Section 7. Let us now focus on the moderate deviations. In general, the setting of Theorem 1 is more challenging since, for $\theta>\pi/4$, we do not see independence in every step. However, for $\theta<\pi/2$, independence can be recovered and we think of this case in the following. In the case $\theta=\pi/2$, we almost surely never see independence of a step from the previous steps and we deal with this case separately in Section 6. More precisely, for $\theta < \pi/2$, the navigation will occasionally make a step that is independent of the past. Conceiving the in-between steps as one segment, we re-enter a regime of independent segments for which we derive the moderate-deviation result. However, the distribution of the segments is hard to trace and has exponential tails. This is the main reason why we restrict ourselves to moderate deviations and only establish the large deviations under some assumptions on the segment distribution in Section 4.

To make this precise, we iteratively define a sequence of *history sets*. We set $H_0 = \emptyset$ and define

$$H_1 := \mathcal{C}_\theta(U_1) \cap B(o, |U_1|),$$

where we recall that our navigation starts at the origin $o$. Here, $B(x, r)$ denotes the open ball with radius $r$ centered at $x \in \mathbb{R}^2$. In words, the history set at step 1 is given by the region that lies both in the (potential) future of the first waypoint $V_1 = U_1$ and in the void space that is responsible for finding the first waypoint. Note that $H_1 = \emptyset$ whenever $0 \le \theta \le \pi/4$. For all larger $n$, we define

$$H_n := \mathcal{C}_\theta(V_{n-1} + U_n) \cap \{H_{n-1} \cup B(V_{n-1}, |U_n|)\},$$

the region in the future of the $n$th waypoint that intersects the joint history. One way to look at the history set $H_n$ is to realize that this is the region where, for the next step, we already know that a certain part of space is already empty of points since we already searched for points there in the previous steps.

Based on our history sets, for $0 < \theta < \pi/2$, we define $\tau_0^\theta = 0$ and

$$\tau_1^\theta := \inf\{n > 0 \colon H_n = \emptyset\},$$

to be the first step where there is no history. For $k > 1$,

$$\tau_k^\theta := \inf\{n > \tau_{k-1}^\theta \colon H_n = \emptyset\}$$

is the $k$th step where the history set is empty. Note that the inter-stopping time gaps $\{\tau_k^\theta - \tau_{k-1}^\theta\}_{k\ge 1}$ are i.i.d. random variables due to the total independence of the underlying Poisson point process. In particular, $\tau_k^\theta = k$ a.s., whenever $0 < \theta \le \pi/4$. We will later verify that $\tau_1^\theta$ is almost-surely finite and has exponential tails. As anticipated, we now build segments

$$U'_i := \sum_{j=\tau_{i-1}^\theta+1}^{\tau_i^\theta} U_j$$

and note that the sequence $\{U'_i\}_{i\ge 1}$ is i.i.d., again due to the total independence of the underlying Poisson point process. Furthermore, let us denote the number of steps before hitting time $t$ by

$$K_t := \sup\left\{n > 0 \colon \sum_{i=1}^{n} X_i < t\right\},$$

where $U_i =: (X_i, Y_i)$ in Cartesian coordinates. Based on our hitting times, we now define

$$K'_t := \sup\{n > 0 \colon \tau_n^\theta \le K_t\} \tag{2}$$

as the index of the largest stopping time before $K_t$. In particular, writing

$$\mathcal{Y}'_t := \sum_{i=1}^{K'_t} Y'_i,$$

where $U'_i =: (X'_i, Y'_i)$ in Cartesian coordinates, we have

$$\sum_{i=1}^{K_t} Y_i = \mathcal{Y}'_t + \sum_{i=\tau^\theta_{K'_t}+1}^{K_t} Y_i. \tag{3}$$

Let us start by noting that $U'_1$ possesses some exponential moments.

**Lemma 2.** *For some* $\gamma = (\gamma_1, \gamma_2)$ *with* $\gamma_1, \gamma_2 > 0$*, we have* $\mathbb{E}[\exp(\langle \gamma, U'_1 \rangle)] < \infty$.

We present the proof of this lemma, along with the proofs of all other statements in this section, later in Section 5. In Section 4, roughly speaking, we make the assumption that $U'_1$ obeys a large-deviation principle, which is the critical input in order to establish the large-deviation principle for $\mathcal{Y}_t$ with $\pi/4 < \theta < \pi/2$.

For the moderate deviations, the difference between $\mathcal{Y}_t$ and $\mathcal{Y}'_t$ is irrelevant in the following sense.

**Proposition 1.** *For any* $0 < \lambda$, $0 < \theta < \pi/2$, *and* $0 < \varepsilon < 1/2$, $\{t^{-1/2-\varepsilon}\mathcal{Y}_t\}_{t\ge 0}$ *and* $\{t^{-1/2-\varepsilon}\mathcal{Y}'_t\}_{t\ge 0}$ *are exponentially equivalent, i.e. for any* $\delta > 0$,

$$\limsup_{t\uparrow\infty} t^{-2\varepsilon} \log \mathbb{P}(|\mathcal{Y}_t - \mathcal{Y}'_t| \ge \delta t^{1/2+\varepsilon}) = -\infty. \tag{4}$$

Now, we want to deal with the randomness in the number of steps $K'_t$ in the definition of $\mathcal{Y}'_t$. For this, we define, for $w > 0$, the process

$$\mathcal{Y}^w_t := \sum_{i=1}^{\lfloor tw \rfloor} Y'_i$$

and establish in the following statement that it obeys the moderate-deviation principle.

**Proposition 2.** *For any* $0 < \lambda$, $0 < \theta < \pi/2$, $0 < w$, *and* $0 < \varepsilon < 1/2$, $\{t^{-1/2-\varepsilon}\mathcal{Y}^w_t\}_{t\ge 0}$ *obeys the moderate-deviation principle with rate function* $I^w_{\lambda,\theta}(x) := x^2/(2w\mathbb{E}[Y'^2_1])$.

Now let $\kappa := \mathbb{E}[X'_1]^{-1}$ denote the inverse expected horizontal progress. By Lemma 2, we have $\kappa > 0$ and, by Lemma 8 below, also $\kappa < \infty$. The following result now presents the final ingredient for the proof of Theorem 1 in the case $0 < \theta < \pi/2$.

**Proposition 3.** *For any* $0 < \lambda$, $0 < \theta < \pi/2$, *and* $0 < \varepsilon < 1/2$, $\{t^{-1/2-\varepsilon}\mathcal{Y}^\kappa_t\}_{t\ge 0}$ *and* $\{t^{-1/2-\varepsilon}\mathcal{Y}'_t\}_{t\ge 0}$ *are exponentially equivalent, i.e. for any* $\delta > 0$,

$$\limsup_{t\uparrow\infty} t^{-2\varepsilon} \log \mathbb{P}(|\mathcal{Y}^\kappa_t - \mathcal{Y}'_t| \ge \delta t^{1/2+\varepsilon}) = -\infty. \tag{5}$$

*Proof of Theorem* 1, *case* $0 < \theta < \pi/2$. By Propositions 1 and 3, we have the exponential equivalence of $\mathcal{Y}_t$ and $\mathcal{Y}^\kappa_t$ and, by Proposition 2, we have that $\mathcal{Y}^\kappa_t$ obeys the moderate-deviation principle with rate function

$$I_{\lambda,\theta}(x) := x^2 \frac{\mathbb{E}[X'_1]}{2\mathbb{E}[{Y'_1}^2]}, \tag{6}$$

as desired. For the scale invariance, if we denote by $U'_{\lambda,1} = (X'_{\lambda,1}, Y'_{\lambda,1})$ the progress of the first independent segment in the navigation based on a Poisson point process with intensity $\lambda > 0$, then, by scaling both the coordinates by $\sqrt{\lambda}$, we find that $(\sqrt{\lambda}X'_{\lambda,1}, \sqrt{\lambda}Y'_{\lambda,1})$ and $(X'_{1,1}, Y'_{1,1})$ have the same distribution. Therefore, we obtain

$$2\rho(\lambda, \theta) = \frac{\mathbb{E}[X'_{\lambda,1}]}{\mathbb{E}[{Y'_{\lambda,1}}^2]} = \frac{\sqrt{\lambda}\mathbb{E}[X'_{1,1}]}{\mathbb{E}[{Y'_{1,1}}^2]} = 2\sqrt{\lambda}\rho(1, \theta),$$

which completes the proof.

Note that, when $0 < \theta \leq \pi/4$, we have $U'_1 = U_1$ and the expression (1) for $\rho$ follows from a simple computation; see, for example, the proof of Lemma 1 in Section 7. The case $\theta = \pi/2$ requires a refined construction of a renewal process. We present the details in Section 6.

## 4. Large Deviations in the Dependent Case

In this section we study the large-deviation behavior of the model in the dependent case, where $\pi/4 < \theta < \pi/2$. Let us start with the large deviations for the empirical average.

**Lemma 3.** *For any $0 < \lambda$ and $\pi/4 < \theta < \pi/2$, $\{n^{-1}\sum_{i=1}^n U'_i\}_{n\geq 1}$ satisfies the large-deviation principle with rate n and rate function*

$$\mathcal{J}'_{\lambda,\theta}(u) := \sup\{\langle\gamma, u\rangle - J'_{\lambda,\theta}(\gamma) \colon \gamma \in \mathbb{R}^2\},$$

*where, for all $\gamma \in \mathbb{R}^2$, $J'_{\lambda,\theta}(\gamma) := \log \mathbb{E}[\exp(\langle\gamma, U'_1\rangle)]$.*

We present the proof of this lemma and the following statements at the end of the paper in Section 7. One would like to establish now a large-deviation principle for $\mathcal{Y}'_t$ similar to the one provided in Theorem 2. However, the situation is more complicated. The key reason for this is that the upper tails of $|U_1|$ are of exponential order $O(-t^2)$ (this is crucially used in step 2 in the proof of Theorem 2), however, the upper tails of $|U'_1|$ are only of exponential order $O(-t)$, as can be seen from the following statement.

**Lemma 4.** *Let $\pi/4 < \theta < \pi/2$, then, for all $n \geq 0$, $\mathbb{P}(\tau_1^\theta > n) \geq ((4\theta - \pi)/(4\theta))^n$.*

We present the proof further below. As a consequence, for example, for $X'_1$, we have, for all $n$ and $s > 0$,

$$\begin{aligned}\mathbb{P}(X'_1 > t) &\geq \mathbb{P}\left(\sum_{i=1}^{\tau_1^\theta} X_i > t,\, \tau_1^\theta > n\right)\\ &\geq \mathbb{P}(\tau_1^\theta > n) - \mathbb{P}\left(\sum_{i=1}^{n} \underline{R}_i \cos\theta \leq t\right)\\ &\geq \mathrm{e}^{an} - \mathrm{e}^{st + n\log\mathbb{E}[\exp(-s\underline{R}\cos\theta)]},\end{aligned}$$

where $a := \log((4\theta - \pi)/(4\theta))$ and $\underline{R}$ and $\{\underline{R}_i\}_{i\geq 1}$ are as in Lemma 5. Hence, we have at least exponential decay if, for $n$ coupled to $t$ as $n = bt$ with $b > 0$, we have $s/b + \log\mathbb{E}[\exp(-s\underline{R}\cos\theta)] < a$. But this is the case by first picking $s$ large and then picking $b$ large. We can proceed similarly for $Y'_1$.

Hence, it is reasonable to assume that we have the following tail behavior within a segment.

**Assumption 1.** *Let* $0 < \lambda$ *and* $\pi/4 < \theta < \pi/2$*. Then, for all* $a > 0$, $b \in \mathbb{R}$, *we have*

$$\lim_{t\uparrow\infty} t^{-1} \log \mathbb{P}\left( \sum_{i=1}^{\lfloor t \rfloor} X_i < at,\ \sum_{i=1}^{\lfloor t \rfloor+1} X_i \geq at,\ \sum_{i=1}^{\lfloor t \rfloor} Y_i \geq bt,\ \tau_1^\theta > t \right) = -\mathcal{H}_{\lambda,\theta}(a, b) \in (-\infty, 0).$$

Note that, it is clear that the rate function must obey our usual scaling relation $\mathcal{H}_{\lambda,\theta} = \sqrt{\lambda}\mathcal{H}_{1,\theta}$. With this in mind, we can prove the large-deviation principle for $\mathcal{Y}_t$.

**Theorem 3.** *Let Assumption 1 hold. Then, for any* $0 < \lambda$ *and* $\pi/4 < \theta < \pi/2$, $\{t^{-1}\mathcal{Y}_t\}_{t\geq 0}$ *obeys the large-deviation principle with rate t and rate function*

$$\mathcal{I}'_{\lambda,\theta}(x) := \inf_{\substack{b\in\mathbb{R},\\ c\in(0,1)}} \big\{ \inf\{\beta \mathcal{J}'_{\lambda,\theta}(c/\beta, b/\beta) \colon \beta > 0\} + \inf\{d\mathcal{H}_{\lambda,\theta}((1-c)/d, (a-b)/d) \colon d > 0\}\big\}.$$

*Moreover,* $\mathcal{I}'_{\lambda,\theta}(x)$ *satisfies the scaling relation* $\mathcal{I}'_{\lambda,\theta}(x) = \sqrt{\lambda}\mathcal{I}'_{1,\theta}(x)$.

Let us try to explain the rate function in words. In the dependent case, the rate function is given by an optimization between the unlikely vertical displacement (up to level $b$, represented by the term involving $\mathcal{J}'_{\lambda,\theta}$), which is cheaper to achieve if the navigation performs fewer steps (up to level $c < 1$) and the, then necessary, cost produced by having an unlikely large last horizontal step (of length of order $(1-c)t$ that covers the remaining vertical displacement, represented by the term involving $\mathcal{H}_{\lambda,\theta}$). The proof is presented in Section 7.

## 5. Proofs for the Moderate Deviations in Case $\theta < \pi/2$

### 5.1. Exponential tails of the stopping times

Let us start by establishing exponential-decay properties of the inter-stopping time gaps.

**Proposition 4.** *Let* $0 < \lambda$ *and* $0 < \theta < \pi/2$*. Then, for all* $n, k \in \mathbb{N}$, *there exist constants* $c, C > 0$ *depending only on* $\theta$ *such that*

$$\mathbb{P}(\tau_k^\theta - \tau_{k-1}^\theta \geq n) \leq C\mathrm{e}^{-cn}. \tag{7}$$

*In particular, for all* $n, k \in \mathbb{N}$, $\mathbb{P}(\tau_k^\theta \geq n) \leq kC\mathrm{e}^{-cn/k}$.

First note that it suffices to prove that $\mathbb{P}(\tau_1^\theta \geq n) \leq C\mathrm{e}^{-cn}$ since the remaining statements then follow from the fact that $\{\tau_k^\theta - \tau_{k-1}^\theta\}_{k\geq 1}$ are i.i.d. We write $\tau := \tau_1^\theta$ for the rest of the proof. It is worth noting that $\tau$ is invariant under scaling both the coordinates by $\sqrt{\lambda}$, and therefore, $\mathbb{P}(\tau \geq n)$ does not depend on $\lambda$. Note also that we can focus on the case $\pi/4 < \theta < \pi/2$ since, for $\theta \leq \pi/4$, we have almost surely that $\tau = 1$ as mentioned before; see the first paragraph in the proof of Theorem 2 in Section 7.

Our arguments rest on ideas developed in [6], more precisely the proof of [6, Proposition 3.1], where however a slightly different situation is analyzed. Let us define, for all $n \geq 0$,

$$L_n := \sup\left\{ x - \sum_{i=1}^{n} X_i \colon (x, y) \in H_n \right\} \vee 0,$$

the *width of the history set* and write $U_n = (X_n, Y_n) = (R_n, \Phi_n)$ for the progress variables first in Cartesian and alternatively in polar coordinates. In particular, $L_0 = 0$ and $\mathbb{P}(\tau \geq n) = \mathbb{P}(H_m \neq \emptyset$ for all $m < n)$.

Since $H_n \subseteq H_{n-1} \cup B(V_{n-1}, R_n)$, we have

$$L_n \le \sup\left\{x - \sum_{i=1}^{n} X_i \colon (x, y) \in H_{n-1} \cup B(V_{n-1}, R_n)\right\} \vee 0$$
$$= \sup\left\{x - \sum_{i=1}^{n} X_i \colon (x, y) \in H_{n-1}\right\} \vee \sup\left\{x - \sum_{i=1}^{n} X_i \colon (x, y) \in B(V_{n-1}, R_n)\right\} \vee 0$$
$$= (L_{n-1} - X_n) \vee (R_n - X_n) \vee 0$$
$$= \max\{(L_{n-1} - X_n)_+, R_n - X_n\},$$

where $x \vee y := \max\{x, y\}$ and $x_+ := \max\{x, 0\}$. Note that, if $\Phi_n \in [-\pi/2 + \theta, \pi/2 - \theta] =: T_\theta$ then the $n$th step escapes its immediate history, i.e. $H_n = \mathcal{C}_\theta(V_n) \cap H_{n-1}$. In that case $H_n \subseteq H_{n-1}$ and, thus, we have

$$L_n \le (L_{n-1} - X_n)_+.$$

Therefore, for all $n \ge 1$, we have

$$L_n \le (L_{n-1} - X_n)_+ \mathbb{1}\{\Phi_n \in T_\theta\} + \max\{(L_{n-1} - X_n)_+, (R_n - X_n)\} \mathbb{1}\{\Phi_n \notin T_\theta\}.$$

The key ingredient for the proof is the following (random) monotone coupling of the progress variables with respect to an i.i.d. sequence, which is also used on multiple occasions throughout the proof section. Denoting $x \wedge y := \min\{x, y\}$, we have the following statement.

**Lemma 5.** *For any $0 < \theta < \pi/2$, there exists a sequence of i.i.d. quadruplets of random variables $\{(\underline{R}_n, \underline{\Phi}_n, \overline{R}_n, \overline{\Phi}_n)\}_{n\ge1}$, defined on an extended probability space, such that, for all $n \ge 1$,*

(i) *$\underline{R}_n \le R_n \le \overline{R}_n$ a.s.,*

(ii) *if $\underline{\Phi}_n \in T_\theta$ then $\Phi_n = \underline{\Phi}_n = \overline{\Phi}_n$ a.s., and*

(iii) *$(\underline{R}_n, \underline{\Phi}_n, \overline{R}_n, \overline{\Phi}_n) \overset{\mathrm{d}}{=} (\underline{R}, \underline{\Phi}, \overline{R}, \overline{\Phi})$, where $(\underline{R}, \underline{\Phi}) := \operatorname{argmin}\{|v| \colon v \in \mathcal{C}_\theta \cap \mathcal{P}_\lambda\}$ and $(\overline{R}, \overline{\Phi}) := \operatorname{argmin}\{|v| \colon v \in \mathcal{C}_{\theta\wedge(\pi/2-\theta)} \cap \mathcal{P}_\lambda\}$, $\underline{\Phi} \sim \mathrm{Uniform}(-\theta, \theta)$ and $\overline{\Phi} \sim \mathrm{Uniform}(-\theta \wedge (\pi/2 - \theta), \theta \wedge (\pi/2 - \theta))$, $\mathbb{P}(\underline{R} > t) = \mathrm{e}^{-\lambda\theta t^2}$ and $\mathbb{P}(\overline{R} > t) = \mathrm{e}^{-\lambda(\theta\wedge(\pi/2-\theta))t^2}$.*

We present the proof of this lemma later in this section. Now, observe that, almost surely,

$$L_n \le (L_{n-1} - R_n \cos\theta)_+ \mathbb{1}\{\Phi_n \in T_\theta\} + \max\{(L_{n-1} - R_n \cos\theta)_+, R_n(1 - \cos\theta)\} \mathbb{1}\{\Phi_n \notin T_\theta\}$$
$$\le (L_{n-1} - \underline{R}_n \cos\theta)_+ \mathbb{1}\{\Phi_n \in T_\theta\} + \max\{L_{n-1}, \overline{R}_n\} \mathbb{1}\{\Phi_n \notin T_\theta\}$$
$$= (L_{n-1} - \underline{R}_n \cos\theta)_+ \mathbb{1}\{\Phi_n \in T_\theta, \underline{\Phi}_n \in T_\theta\} + (L_{n-1} - \underline{R}_n \cos\theta)_+ \mathbb{1}\{\Phi_n \in T_\theta, \underline{\Phi}_n \notin T_\theta\}$$
$$+ \max\{L_{n-1}, \overline{R}_n\} \mathbb{1}\{\Phi_n \notin T_\theta, \underline{\Phi}_n \in T_\theta\} + \max\{L_{n-1}, \overline{R}_n\} \mathbb{1}\{\Phi_n \notin T_\theta, \underline{\Phi}_n \notin T_\theta\}$$
$$\le (L_{n-1} - \underline{R}_n \cos\theta)_+ \mathbb{1}\{\underline{\Phi}_n \in T_\theta\} + \max\{L_{n-1}, \overline{R}_n\} \mathbb{1}\{\underline{\Phi}_n \notin T_\theta\}.$$

Note that, by Lemma 5(ii), we have $\{\underline{\Phi}_n \in T_\theta\} \subseteq \{\Phi_n \in T_\theta\}$, which, together with the fact that $L_{n-1} - \underline{R}_n \cos\theta \le L_{n-1} \le \max\{L_{n-1}, \overline{R}_n\}$, implies the final inequality. We now make a comparison with a Markov chain $\{M_n\}_{n\geq 0}$ on $\mathbb{N}_0$, defined as $M_0 := 0$ and

$$M_n := (M_{n-1} - \lfloor \underline{R}_n \cos\theta \rfloor)_+ \mathbb{1}\{\underline{\Phi}_n \in T_\theta\} + \max\{M_{n-1}, \lceil \overline{R}_n \rceil\} \mathbb{1}\{\underline{\Phi}_n \notin T_\theta\}, \quad n > 0.$$

The following result establishes that the stopping time $\tau^M := \inf\{n > 0 \colon M_n = 0\}$ of the Markov chain to hit 0 has an exponential tail.

**Lemma 6** *For any $n > 0$, $\mathbb{P}(\tau^M \geq n) \leq C\mathrm{e}^{-cn}$ for some $C, c > 0$.*

*Proof of Proposition* 4. The reason for constructing the Markov chain is to dominate the sequence of widths $\{L_n\}_{n\geq 0}$. More precisely, by construction, $M_n \geq L_n$ for all $n \geq 0$ and, in particular, when $M_n = 0$, we have $L_n = 0$. Therefore, almost surely, $\tau^M \geq \tau$, which completes the proof.

*Proof of Lemma* 6. Note that the Markov chain $\{M_n\}_{n\geq 0}$ is irreducible since, for any $m > 0$,

$$\begin{aligned} &\mathbb{P}(M_1 = m \mid M_0 = 0) \geq \mathbb{P}(\lceil \overline{R}_1 \rceil = m \text{ and } \underline{\Phi}_1 \notin T_\theta) > 0 \quad \text{and} \\ &\mathbb{P}(M_1 = 0 \mid M_0 = m) \geq \mathbb{P}(\lfloor \underline{R}_1 \cos\theta \rfloor > m \text{ and } \underline{\Phi}_1 \in T_\theta) > 0. \end{aligned} \tag{8}$$

Now, we demonstrate the recurrence of $\{M_n\}_{n\geq 0}$ by establishing that 0 is a recurrent state. From Lemma 5, we know that, for all $h \geq 2$,

$$\mathbb{P}(\lceil \overline{R}_1 \rceil > h) \leq \mathbb{P}(\overline{R}_1 > h/2) = \exp(-\lambda(\pi/2 - \theta)h^2/4).$$

Let us define the sequence $\{a_n\}_{n\geq 1}$ as

$$a_n := \left\lceil \left(\frac{8}{\lambda(\pi/2 - \theta)} \log n\right)^{1/2} \right\rceil.$$

Then, for all sufficiently large $n$, we have $\mathbb{P}(\lceil \overline{R}_1 \rceil > a_n) \leq n^{-2}$. This implies that, for all sufficiently large $n$,

$$\mathbb{P}\Big(\max_{i=1}^{n} \lceil \overline{R}_i \rceil \leq a_n\Big) \geq (1 - n^{-2})^n$$

and, therefore,

$$\lim_{n\uparrow\infty} \mathbb{P}\Big(\max_{i=1}^{n} \lceil \overline{R}_i \rceil \leq a_n\Big) = 1. \tag{9}$$

Now, let $G_n := \mathbb{1}\{\underline{\Phi}_n \in T_\theta \text{ and } \underline{R}_n \cos\theta \geq 1\}$ and $q := \mathbb{P}(G_1 = 1) > 0$. From the definition of $M_n$, it follows that, whenever $G_n = 1$, we have $M_n \leq (M_{n-1} - 1)_+$. We define $\mathcal{A}_n$ to be the event that the finite sequence $\{G_i\}_{i=1}^n$ has a run of 1s of length at least $a_n$. Since $\{G_i\}_{i=1}^n$ are i.i.d., we have

$$\mathbb{P}(\mathcal{A}_n^c) \leq (1 - q^{a_n})^{\lfloor n/a_n \rfloor}.$$

Furthermore, since, for all sufficiently large $n$,

$$a_n \leq \frac{1}{2(-\log q)} \log n,$$

we have, for all sufficiently large $n$,

$$\mathbb{P}(\mathcal{A}_n^c) \le (1 - n^{-1/2})^{\lfloor 2(-\log q)n/\log n \rfloor},$$

which tends to 0 as $n \to \infty$. Therefore, we obtain

$$\lim_{n\uparrow\infty} \mathbb{P}(\mathcal{A}_n) = 1. \tag{10}$$

Now, by construction, we have

$$\begin{aligned}
&\mathbb{P}(M_i = 0 \text{ for some } 1 \le i \le n \mid M_0 = 0) \\
&\quad \ge \mathbb{P}\Big(\max_{i=1}^{n} M_i \le a_n \text{ and } \{G_i\}_{i=1}^n \text{ has } a_n \text{ consecutive } 1\text{ s} \;\Big|\; M_0 = 0\Big) \\
&\quad \ge \mathbb{P}\Big(\max_{i=1}^{n} \lceil \overline{R}_i \rceil \le a_n \text{ and } \{G_i\}_{i=1}^n \text{ has } a_n \text{ consecutive } 1\text{ s} \;\Big|\; M_0 = 0\Big) \\
&\quad = \mathbb{P}\Big(\max_{i=1}^{n} \lceil \overline{R}_i \rceil \le a_n \text{ and } \{G_i\}_{i=1}^n \text{ has } a_n \text{ consecutive } 1\text{ s}\Big).
\end{aligned}$$

But this tends to 1 as $n \to \infty$ by (9) and (10). As a result, we established that 0 is a recurrent state, and therefore, in view of (8), the Markov chain $\{M_n\}_{n\ge 0}$ is recurrent.

The rest of the proof is similar to the proof of [2, Proposition 5.5], but we include some details for the convenience of the reader. Observe that

$$\begin{aligned}
&\mathbb{E}[\mathrm{e}^{M_1} \mid M_0 = k] \\
&\quad = \mathbb{E}\big[\exp\big((k - \lfloor \underline{R}_1 \cos\theta \rfloor)_+ \mathbb{1}\{\underline{\Phi}_1 \in T_\theta\} + \max\{k, \lceil \overline{R}_1 \rceil\} \mathbb{1}\{\underline{\Phi}_1 \notin T_\theta\}\big)\big] \qquad (11) \\
&\quad = \mathrm{e}^k \mathbb{E}\big[\exp\big(-\min\{k, \lfloor \underline{R}_1 \cos\theta \rfloor\} \mathbb{1}\{\underline{\Phi}_1 \in T_\theta\} + \max\{0, \lceil \overline{R}_1 \rceil - k\} \mathbb{1}\{\underline{\Phi}_1 \notin T_\theta\}\big)\big]
\end{aligned}$$

and note that

$$\exp\big(-\min\{k, \lfloor \underline{R}_1 \cos\theta \rfloor\} \mathbb{1}\{\underline{\Phi}_1 \in T_\theta\} + \max\{0, \lceil \overline{R}_1 \rceil - k\} \mathbb{1}\{\underline{\Phi}_1 \notin T_\theta\}\big) \le \exp(\lceil \overline{R}_1 \rceil),$$

which has finite expectation. Therefore, for all $k \ge 0$,

$$\mathbb{E}[\exp(M_1) \mid M_0 = k] < \infty. \tag{12}$$

Moreover, by the dominated-convergence theorem, we get

$$\begin{aligned}
&\lim_{k\uparrow\infty} \mathbb{E}[\exp(-\min\{k, \lfloor \underline{R}_1 \cos\theta \rfloor\} \mathbb{1}\{\underline{\Phi}_1 \in T_\theta\} + \max\{0, \lceil \overline{R}_1 \rceil - k)\} \mathbb{1}\{\underline{\Phi}_1 \notin T_\theta\})] \\
&\qquad = \mathbb{E}[\exp(-\lfloor \underline{R}_1 \cos\theta \rfloor \mathbb{1}\{\underline{\Phi}_1 \in T_\theta\})] \\
&\qquad < 1
\end{aligned}$$

and, thus, there exists $r > 1$ and $k_0 \ge 0$ such that, for all $k > k_0$,

$$\mathbb{E}[\exp(-\min\{k, \lfloor \underline{R}_1 \cos\theta \rfloor\} \mathbb{1}\{\underline{\Phi}_1 \in T_\theta\} + \max\{0, \lceil \overline{R}_1 \rceil - k)\} \mathbb{1}\{\underline{\Phi}_1 \notin T_\theta\})] < r^{-1}.$$

But this, together with (11), implies that, for all $k > k_0$,

$$\mathbb{E}[\exp(M_1) \mid M_0 = k] < r^{-1} \exp(k). \tag{13}$$

Let $\sigma := \inf\{n > 0 \colon M_n \le k_0\}$ denote the first time that $M_n$ is below $k_0$ and fix any $l > k_0$. In view of (13), note that, for $M_0 = l$, the process $\{Y_n\}_{n\ge 0}$ defined as $Y_n := r^{n\wedge\sigma} \exp(M_{n\wedge\sigma})$ is a nonnegative supermartingale. Therefore, by the recurrence of the Markov chain $\{M_n\}_{n\ge 0}$, almost surely,

$$Y_n \to Y_\infty := r^\sigma \exp(M_\sigma) \ge r^\sigma.$$

This, together with Fatou's lemma, implies that

$$\mathbb{E}[r^\sigma \mid M_0 = l] \le \mathbb{E}[Y_\infty \mid M_0 = l] \le \liminf_{n\uparrow\infty} \mathbb{E}[Y_n \mid M_0 = l] \le \mathbb{E}[Y_0 \mid M_0 = l] = \exp(l).$$

Since $l > k_0$ is arbitrary, we have, for all $l > k_0$,

$$\mathbb{E}[r^\sigma | M_0 = l] \le \exp(l).$$

Therefore, for any $z \le k_0$,

$$\begin{aligned}
\mathbb{E}[r^\sigma \mid M_0 = z] &\le r + r\sum_{l>k_0} \mathbb{P}(M_1 = l \mid M_0 = z)\mathbb{E}[r^\sigma | M_0 = l] \\
&\le r + r\sum_{l>k_0} \mathbb{P}(M_1 = l \mid M_0 = z)\exp(l) \\
&\le r + r\mathbb{E}[\exp(M_1) \mid M_0 = z] \\
&< \infty.
\end{aligned}$$

Note that, for any $z \le k_0$ and $l > k_0$,

$$\mathbb{P}(M_1 = l \mid M_0 = z) = \mathbb{P}(\underline{\Phi}_1 \notin T_\theta, \lceil \overline{R}_1 \rceil = l) = \mathbb{P}(M_1 = l \mid M_0 \le k_0),$$

which does not depend on $z$. Therefore, if we define $\{\sigma_i\}_{i\ge 0}$ as $\sigma_0 := 0$ and, for all $i \ge 1$,

$$\sigma_i := \inf\{n > \sigma_{i-1} \colon M_n \le k_0\},$$

then $\{\sigma_i - \sigma_{i-1}\}_{i\ge 1}$ are i.i.d. copies of $\sigma$. Moreover, for all $z \le k_0$,

$$c_1 := \mathbb{E}[r^\sigma \mid M_0 \le k_0] = \mathbb{E}[r^\sigma \mid M_0 = z] < \infty.$$

We choose $c_2 > 0$ such that ${c_1}^{c_2} < r$. Then, by Markov's inequality,

$$\mathbb{P}(\sigma_{\lfloor c_2 n\rfloor} \ge n \mid M_0 = 0) \le r^{-n}\mathbb{E}[r^{\sigma_{\lfloor c_2 n\rfloor}} \mid M_0 = 0] = r^{-n}\mathbb{E}[r^{\sigma_1} \mid M_0 = 0]^{\lfloor c_2 n\rfloor} \le (r^{-1}{c_1}^{c_2})^n,$$

and thus,

$$\begin{aligned}
&\mathbb{P}(\tau^M > n \mid M_0 = 0) \\
&\quad\le \mathbb{P}(\tau^M > n, \sigma_{\lfloor c_2 n\rfloor} < n \mid M_0 = 0) + \mathbb{P}(\sigma_{\lfloor c_2 n\rfloor} \ge n \mid M_0 = 0) \\
&\quad\le \mathbb{P}\left(\bigcap_{i=1}^{\lfloor c_2 n\rfloor} \{\underline{\Phi}_{\sigma_i+1} \notin T_\theta \text{ or } \lfloor \underline{R}_{\sigma_i+1} \cos\theta \rfloor < k_0\}\right) + \mathbb{P}(\sigma_{\lfloor c_2 n\rfloor} \ge n \mid M_0 = 0) \\
&\quad\le \mathbb{P}(\underline{\Phi}_1 \notin T_\theta \text{ or } \lfloor \underline{R}_1 \cos\theta \rfloor < k_0)^{\lfloor c_2 n\rfloor} + (r^{-1}{c_1}^{c_2})^n \\
&\quad\le C_1 \mathrm{e}^{-cn}
\end{aligned}$$

for some $C_1, c > 0$. This proves the result.

Now, the only thing that remains to be proven is Lemma 5. For this, let us first verify the following statement.

**Lemma 7.** *For any $n \geq 0$ and $\pi/4 \leq \theta < \pi/2$, we have $\mathcal{C}_{\pi/2-\theta}(V_n) \cap H_n = \emptyset$.*

*Proof of Lemma 7.* From the definition of $H_n$, we see that, for any $n \geq 0$,

$$H_n \subseteq \bigcup_{m=0}^{n-1} B(V_m, R_{m+1}).$$

Now, since, for any $0 \leq m \leq n-1$, $V_n \in \mathcal{C}_\theta(V_m)$, it follows that $R_{m+1} \leq |V_n - V_m|$. Therefore,

$$H_n \subseteq \bigcup_{m=0}^{n-1} B(V_m, |V_n - V_m|).$$

Hence, it suffices to show that, for any $0 \leq m \leq n-1$, $\mathcal{C}_{\pi/2-\theta}(V_n) \cap B(V_m, |V_n - V_m|) = \emptyset$. For this, note that the lines

$$L_n(\pi/2-\theta) := V_n + \{(r, \varphi) \colon r > 0, \varphi = \pi/2 - \theta\} \quad \text{and}$$
$$L_n(-\pi/2+\theta) := V_n + \{(r, \varphi) \colon r > 0, \varphi = -\pi/2 + \theta\}$$

are the two boundary lines of $\mathcal{C}_{\pi/2-\theta}(V_n)$. Hence, it is also adequate to prove that, for any $0 \leq m \leq n-1$, $(L_n(\pi/2-\theta) \cup L_n(-\pi/2+\theta)) \cap B(V_m, |V_n - V_m|) = \emptyset$.

We now fix $m < n$ and write $V_m$ and $V_n$ in Cartesian coordinates as $V_m = (V_{m,1}, V_{m,2})$ and $V_n = (V_{n,1}, V_{n,2})$. If $V_{m,2} = V_{n,2}$ then $|V_n - V_m| = V_{n,1} - V_{m,1}$ and since, for any $x = (x_1, x_2) \in L_n(\pi/2-\theta) \cup L_n(-\pi/2+\theta)$, we have $x_1 > V_{n,1} > V_{m,1}$, it follows that

$$|x - V_m| \geq x_1 - V_{m,1} > V_{n,1} - V_{m,1} = |V_n - V_m|.$$

Therefore, $(L_n(\pi/2-\theta) \cup L_n(-\pi/2+\theta)) \cap B(V_m, |V_n - V_m|) = \emptyset$.

Now, suppose that $V_{m,2} \neq V_{n,2}$. Without loss of generality, we assume that $V_{m,2} < V_{n,2}$. Then, for any $x = (x_1, x_2) \in L_n(\pi/2-\theta)$, we have $x_1 > V_{n,1} > V_{m,1}$ and $x_2 > V_{n,2} > V_{m,2}$, which implies that

$$|x - V_m| = \sqrt{(x_1 - V_{m,1})^2 + (x_2 - V_{m,2})^2} > \sqrt{(V_{n,1} - V_{m,1})^2 + (V_{n,2} - V_{m,2})^2} = |V_n - V_m|.$$

Therefore, $L_n(\pi/2-\theta) \cap B(V_m, |V_n - V_m|) = \emptyset$. Now, we draw two horizontal line segments $\overline{V_m A}$ and $\overline{V_n B}$ passing through $V_m$ and $V_n$, respectively, such that $\overline{V_m A}$ intersects $L_n(-\pi/2+\theta)$ at point $C$; see Figure 2 for an illustration. From the definition of $L_n(-\pi/2+\theta)$, we know that $\angle B V_n C = \pi/2 - \theta$. Since $\overline{V_m A}$ and $\overline{V_n B}$ are parallel, we also have $\angle V_n C V_m = \pi/2 - \theta$. Furthermore, as $V_n \in \mathcal{C}_\theta(V_m)$, it follows that $\angle V_n V_m C \leq \theta$. Now, focusing on the triangle $\triangle V_n V_m C$, we observe the relationship

$$\angle V_n C V_m + \angle V_n V_m C + \angle V_m V_n C = \pi,$$

which implies that $\angle V_m V_n C \geq \pi/2$, and this, in turn, implies that $\cos \angle V_m V_n C \leq 0$. Therefore, for any $x \in L_n(-\pi/2+\theta)$, by the law of cosines,

$$|x - V_m| = \sqrt{|V_n - V_m|^2 + |x - V_n|^2 - 2\cos \angle V_m V_n C \cdot |V_n - V_m| \cdot |x - V_n|} > |V_n - V_m|,$$

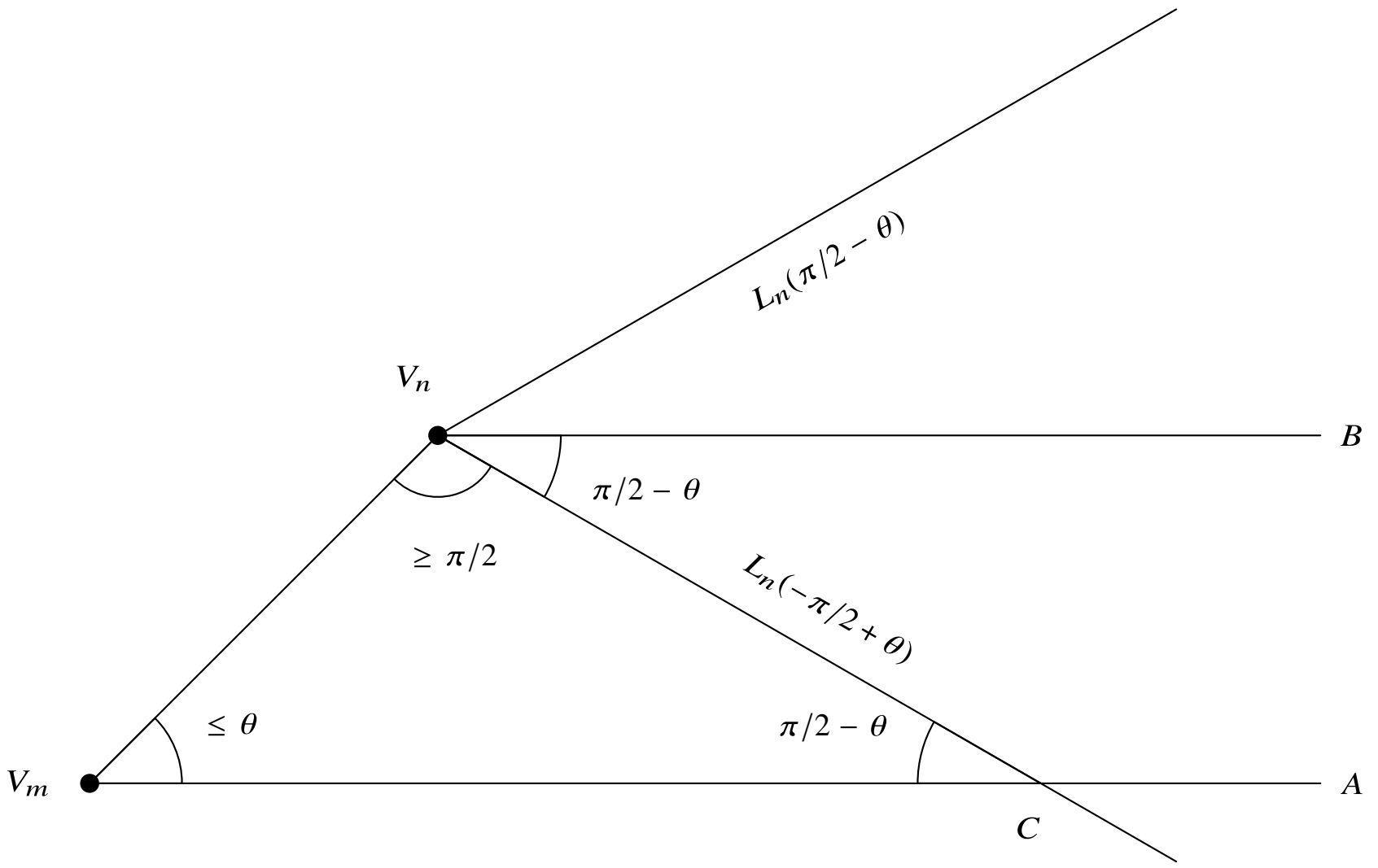

FIGURE 2. Illustration for the proof of Lemma 7.

which means that $L_n(\pi/2-\theta)\cap B(V_m, |V_n - V_m|) = \emptyset$. This proves that, for any $0 \leq m \leq n-1$, $\mathcal{C}_{\pi/2-\theta}(V_n)\cap B(V_m, |V_n - V_m|) = \emptyset$, thereby completing the proof.

*Proof of Lemma* 5. For $0 < \theta \leq \pi/4$, the random variables $\{(R_n, \Phi_n)\}_{n\geq 1}$ are i.i.d. Therefore, Lemma 5 holds trivially by taking $(\underline{R}_n, \underline{\Phi}_n) = (\overline{R}_n, \overline{\Phi}_n) = (R_n, \Phi_n)$ for all $n \geq 1$. So, we now assume that $\pi/4 < \theta < \pi/2$. We define a sequence of sets $\{D_n\}_{n\geq 1}$ as

$$D_n := \mathcal{C}_\theta(V_{n-1}) \cap \overline{B(V_{n-1}, R_n)} \cap H^c_{n-1}.$$

Essentially, $D_n$ is the previously unexplored set where we have explored at the $n$th step to find $V_n$. Clearly, $\{D_n\}_{n\geq 1}$ are pairwise disjoint. Let $\{\mathcal{P}^{(n)}_\lambda\}_{n\geq 1}$ be i.i.d. copies of the Poisson point process $\mathcal{P}_\lambda$. We define $\{\mathcal{Q}^{(n)}_\lambda\}_{n\geq 1}$ as

$$\mathcal{Q}^{(n)}_\lambda := (\mathcal{P}_\lambda \cap D_n) \cup (\mathcal{P}^{(n)}_\lambda \cap D^c_n).$$

Since the $D_n$ are disjoint, $\{\mathcal{Q}^{(n)}_\lambda\}_{n\geq 1}$ is a sequence of i.i.d. Poisson point processes with intensity $\lambda$. We write

$$\underline{W}_n := \operatorname{argmin}\{|v - V_{n-1}| : v \in \mathcal{C}_\theta(V_{n-1}) \cap \mathcal{Q}^{(n)}_\lambda\},$$

and define $(\underline{R}_n, \underline{\Phi}_n)$ to be the polar coordinates of $\underline{W}_n - V_{n-1}$. Similarly, we write

$$\overline{W}_n := \operatorname{argmin}\{|v - V_{n-1}| : v \in \mathcal{C}_{\pi/2-\theta}(V_{n-1}) \cap \mathcal{Q}^{(n)}_\lambda\},$$

and define $(\overline{R}_n, \overline{\Phi}_n)$ to be the polar coordinates of $\overline{W}_n - V_{n-1}$. From the definition, it follows that $(\underline{R}_n, \underline{\Phi}_n, \overline{R}_n, \overline{\Phi}_n)$ and $(\underline{R}, \underline{\Phi}, \overline{R}, \overline{\Phi})$ are equal in distribution. Furthermore, note that, since $\{\mathcal{Q}^{(n)}_\lambda\}_{n\geq 1}$ are i.i.d., $\{(\underline{R}_n, \underline{\Phi}_n, \overline{R}_n, \overline{\Phi}_n)\}_{n\geq 1}$ are also i.i.d. We know that $(R_n, \Phi_n)$ are the polar coordinates of $V_n - V_{n-1}$, where

$$\begin{aligned}
V_n &= \operatorname{argmin}\{|v - V_{n-1}| \colon v \in \mathcal{C}_\theta(V_{n-1}) \cap \mathcal{P}_\lambda\} \\
&= \operatorname{argmin}\{|v - V_{n-1}| \colon v \in \mathcal{C}_\theta(V_{n-1}) \cap \overline{B(V_{n-1}, R_n)} \cap H^c_{n-1} \cap \mathcal{P}_\lambda\} \\
&= \operatorname{argmin}\{|v - V_{n-1}| \colon v \in D_n \cap \mathcal{P}_\lambda\} \\
&= \operatorname{argmin}\{|v - V_{n-1}| \colon v \in D_n \cap \mathcal{Q}^{(n)}_\lambda\} \\
&= \operatorname{argmin}\{|v - V_{n-1}| \colon v \in \mathcal{C}_\theta(V_{n-1}) \cap \overline{B(V_{n-1}, R_n)} \cap H^c_{n-1} \cap \mathcal{Q}^{(n)}_\lambda\} \\
&= \operatorname{argmin}\{|v - V_{n-1}| \colon v \in \mathcal{C}_\theta(V_{n-1}) \cap H^c_{n-1} \cap \mathcal{Q}^{(n)}_\lambda\}.
\end{aligned}$$

Now, from Lemma 7, we obtain $\mathcal{C}_{\pi/2-\theta}(V_{n-1}) \subseteq \mathcal{C}_\theta(V_{n-1}) \cap H^c_{n-1} \subseteq \mathcal{C}_\theta(V_{n-1})$, which implies that $\overline{R}_n \ge R_n \ge \underline{R}_n$. Additionally, if $\underline{\Phi}_n \in T_\theta$ then we have $(\underline{R}_n, \underline{\Phi}_n) = (R_n, \Phi_n) = (\overline{R}_n, \overline{\Phi}_n)$. This concludes the proof.

### 5.2. Proof of Lemma 2 and Proposition 2

*Proof of Lemma* 2. Using Proposition 4, Lemma 5, and the Cauchy–Schwarz inequality, we get

$$\mathbb{E}[\mathrm{e}^{\langle \gamma, U'_1\rangle}] \le \sum_{n\ge 1} \mathbb{E}[\mathrm{e}^{2\langle \gamma, \sum_{j=1}^n U_j\rangle}]^{1/2} \mathbb{P}(\tau^\theta_1 \ge n)^{1/2} \le \sum_{n\ge 1} \mathbb{E}[\mathrm{e}^{2(\gamma_1+\gamma_2)\overline{R}}]^{n/2} \sqrt{C} \mathrm{e}^{-cn/2}.$$

Since $\overline{R}$ has all exponential moments finite, for all sufficiently small $\gamma_1, \gamma_2 > 0$, the above sum is also finite.

*Proof of Proposition* 2. By symmetry and the exponential Markov inequality, for sufficiently small $s > 0$,

$$\mathbb{P}(|Y'_1| \ge n^{1/2+\varepsilon}) \le 2\mathrm{e}^{-sn^{1/2+\varepsilon}} \mathbb{E}[\mathrm{e}^{sY'_1}] < \infty,$$

where we also used Lemma 2, and thus,

$$\limsup_{n\uparrow\infty} n^{-2\varepsilon} \log \mathbb{P}(|Y'_1| \ge n^{1/2+\varepsilon}) \le -s \limsup_{n\uparrow\infty} n^{1/2-\varepsilon} = -\infty.$$

Hence, $\{\mathcal{Y}^w_{\lfloor t\rfloor}\}_{t\ge 0}$ satisfies the conditions of [9, Theorem 2.2] and, thus, obeys the moderate-deviation principle with rate $x^2/(2w\mathbb{E}[Y'^2_1])$. But $\{\mathcal{Y}^w_{\lfloor t\rfloor}\}_{t\ge 0}$ and $\{\mathcal{Y}^w_t\}_{t\ge 0}$ are exponentially equivalent since, for all $\delta > 0$,

$$\limsup_{t\uparrow\infty} t^{-2\varepsilon} \log \mathbb{P}(|\mathcal{Y}^w_{\lfloor t\rfloor} - \mathcal{Y}^w_t| \ge \delta t^{1/2+\varepsilon}) \le \limsup_{t\uparrow\infty} t^{-2\varepsilon} \log \mathbb{P}\Bigg(\sum_{i=1}^{\lceil w+1\rceil} |Y'_i| \ge \delta t^{1/2+\varepsilon}\Bigg),$$

as $\lfloor tw\rfloor - \lfloor\lfloor t\rfloor w\rfloor \le \lceil w+1\rceil$. Hence, again using symmetry, independence, and the exponential Markov inequality, for sufficiently small $s > 0$,

$$\mathbb{P}\Bigg(\sum_{i=1}^{\lceil w+1\rceil} |Y'_i| \ge \delta t^{1/2+\varepsilon}\Bigg) \le \mathrm{e}^{-s\delta t^{1/2+\varepsilon}} \mathbb{E}[\mathrm{e}^{s|Y'_1|}]^{\lceil w+1\rceil} \le \mathrm{e}^{-s\delta t^{1/2+\varepsilon}} 2^{\lceil w+1\rceil} \mathbb{E}[\mathrm{e}^{sY'_1}]^{\lceil w+1\rceil} < \infty,$$

and thus,

$$\limsup_{t\uparrow\infty} t^{-2\varepsilon}\log \mathbb{P}(|\mathcal{Y}^w_{\lfloor t\rfloor} - \mathcal{Y}^w_t| \geq \delta t^{1/2+\varepsilon}) \leq -s\delta \limsup_{t\uparrow\infty} t^{1/2-\varepsilon} = -\infty,$$

as desired.

## 5.3 Proof of Propositions 1

*Proof of Proposition* 1. First note that we have $\mathcal{Y}_t = \sum_{i=1}^{K_t} Y_i + Y_{K_t+1}(t - \sum_{i=1}^{K_t} X_i)/X_{K_t+1}$ with $(t - \sum_{i=1}^{K_t} X_i)/X_{K_t+1} \leq 1$, and hence, by Lemma 5,

$$\begin{aligned}
\mathbb{P}(|\mathcal{Y}_t - \mathcal{Y}'_t| \geq \delta t^{1/2+\varepsilon}) &\leq \mathbb{P}\left(\sum_{i=\tau^\theta_{K'_t}+1}^{K_t+1} |Y_i| \geq \delta t^{1/2+\varepsilon}\right) \\
&\leq \mathbb{P}\left(\sum_{i=\tau^\theta_{K'_t}+1}^{\tau^\theta_{K'_t+1}} |Y_i| \geq \delta t^{1/2+\varepsilon}\right) \\
&\leq \sum_{j=0}^{\lfloor t^2\rfloor} \mathbb{P}\left(\sum_{i=\tau^\theta_j+1}^{\tau^\theta_{j+1}} |Y_i| \geq \delta t^{1/2+\varepsilon},\, K'_t = j\right) + \mathbb{P}(K'_t \geq t^2) \\
&\leq t^2\mathbb{P}\left(\sum_{i=1}^{\tau^\theta_1} |Y_i| \geq \delta t^{1/2+\varepsilon}\right) + \mathbb{P}(K'_t \geq t^2) \\
&\leq t^2\mathbb{P}\left(\sum_{i=1}^{\lfloor t^{\varepsilon'}\rfloor} \overline{R}_i \geq \delta t^{1/2+\varepsilon}\right) + t^2\mathbb{P}(\tau^\theta_1 > \lfloor t^{\varepsilon'}\rfloor) + \mathbb{P}\left(\sum_{i=1}^{\lceil t^2\rceil} X'_i \leq t\right)
\end{aligned}$$

for some $\varepsilon' \in (2\varepsilon, 1/2+\varepsilon)$. This, together with Proposition 4 and the exponential Markov inequality, implies that

$$\mathbb{P}(|\mathcal{Y}_t - \mathcal{Y}'_t| \geq \delta t^{1/2+\varepsilon}) \leq t^2 \mathrm{e}^{-\delta t^{1/2+\varepsilon}} \mathbb{E}[\mathrm{e}^{\overline{R}}]^{\lfloor t^{\varepsilon'}\rfloor} + t^2 C\mathrm{e}^{-c\lfloor t^{\varepsilon'}\rfloor} + \mathrm{e}^t \mathbb{E}[\mathrm{e}^{-X'_1}]^{\lceil t^2\rceil},$$

and therefore,

$$\begin{aligned}
&\limsup_{t\uparrow\infty} t^{-2\varepsilon}\log \mathbb{P}(|\mathcal{Y}_t - \mathcal{Y}'_t| \geq \delta t^{1/2+\varepsilon}) \\
&\quad\leq -\min\Big\{\liminf_{t\uparrow\infty} t^{-2\varepsilon}(\delta t^{1/2+\varepsilon} - \lfloor t^{\varepsilon'}\rfloor \log \mathbb{E}[\mathrm{e}^{\overline{R}}]),\, \liminf_{t\uparrow\infty} c\lfloor t^{\varepsilon'}\rfloor t^{-2\varepsilon}, \\
&\qquad\qquad -\limsup_{t\uparrow\infty} t^{-2\varepsilon}(t + \lceil t^2\rceil \log \mathbb{E}[\mathrm{e}^{-X'_1}])\Big\} \\
&\quad= -\infty,
\end{aligned}$$

as desired.

### 5.4. Proof of Proposition 3

Before we prove Proposition 3, let us collect concentration properties of the horizontal displacement.

**Lemma 8.** *We have* $\kappa := \mathbb{E}[X_1']^{-1} \in (0, \infty)$ *and, for all* $\varepsilon > 0$*, there exist constants* $c_1, c_2 > 0$ *such that, for all* $t > 0$*,*

$$\mathbb{P}(K_t'/t \notin B_\varepsilon(\kappa)) \le c_1 \mathrm{e}^{-c_2 t}. \tag{14}$$

*Proof of Lemma* 8. First, note that

$$\begin{aligned}
\kappa^{-1} &= \mathbb{E}[X_1'] \\
&\ge \mathbb{E}[X_1 \mathbb{1}\{\tau_1^\theta = 1\}] \\
&= \lambda \int_0^\infty \mathrm{d}r\, r^2 \mathrm{e}^{-\theta r^2} \int_{-\theta}^{\theta} \mathrm{d}\varphi\, \cos\varphi \mathbb{1}\left\{\theta - \frac{\pi}{2} \le \varphi \le \frac{\pi}{2} - \theta\right\} \\
&> 0,
\end{aligned}$$

since independence after one step is achieved precisely if the angle is in the interval $[\theta - \pi/2, \pi/2 - \theta]$. Using Lemma 5 and Proposition 4 together with the Cauchy–Schwarz inequality, we get

$$\begin{aligned}
\kappa^{-1} &= \mathbb{E}[X_1'] \\
&= \mathbb{E}\left[\sum_{i=1}^{\tau} X_i\right] \\
&\le \mathbb{E}\left[\sum_{i=1}^{\tau} \overline{R}_i\right] \\
&\le \sum_{i=1}^{\infty} \mathbb{E}[\overline{R}_i \mathbb{1}\{\tau \ge i\}] \\
&\le \sum_{i=1}^{\infty} (\mathbb{E}[\overline{R}_i^2]\mathbb{P}(\tau \ge i))^{1/2} \\
&\le \sum_{i=1}^{\infty} (\mathbb{E}[\overline{R}^2]C)^{1/2} \mathrm{e}^{-ci/2} \\
&< \infty.
\end{aligned}$$

For the second part of the statement, note that we can bound

$$\begin{aligned}
\mathbb{P}(K_t'/t \notin B_\varepsilon(\kappa)) &\le \mathbb{P}(K_t' \le t(\kappa - \varepsilon)) + \mathbb{P}(K_t' \ge t(\kappa + \varepsilon)) \\
&\le \mathbb{P}\left(\sum_{i=1}^{\lceil t(\kappa-\varepsilon)\rceil} X_i' \ge t\right) + \mathbb{P}\left(\sum_{i=1}^{\lfloor t(\kappa+\varepsilon)\rfloor} X_i' < t\right),
\end{aligned} \tag{15}$$

where we used the fact that $K'_t$, as defined in (2), can also be represented as

$$K'_t = \sup\left\{n > 0\colon \sum_{i=1}^{n} X'_i < t\right\}.$$

Now, by the exponential Markov inequality, for any $s > 0$,

$$\limsup_{t\uparrow\infty} t^{-1} \log \mathbb{P}\left(\sum_{i=1}^{\lceil t(\kappa-\varepsilon)\rceil} X'_i \geq t\right) \leq -s(1-(\kappa-\varepsilon)s^{-1}\log \mathbb{E}[\mathrm{e}^{sX'_1}]).$$

Using Lemma 5 and Proposition 4 together with the Cauchy–Schwarz inequality, we obtain

$$\begin{aligned}
\mathbb{E}[\mathrm{e}^{sX'_1}] &\leq \mathbb{E}[\mathrm{e}^{\sum_{i=1}^{\tau} s\bar{R}_i}] \\
&= \sum_{n=1}^{\infty} \mathbb{E}[\mathrm{e}^{\sum_{i=1}^{n} s\bar{R}_i} \mathbb{1}\{\tau = n\}] \\
&\leq \sum_{n=1}^{\infty} (\mathbb{E}[\mathrm{e}^{\sum_{i=1}^{n} 2s\bar{R}_i}]\mathbb{P}(\tau = n))^{1/2} \\
&\leq \sum_{n=1}^{\infty} (\mathbb{E}[\mathrm{e}^{2s\bar{R}}])^{n/2}\sqrt{C}\mathrm{e}^{-cn/2},
\end{aligned}$$

which is finite for all $0 < s \leq s_0$ and for some $s_0 > 0$, as a consequence of Lemma 5. Hence, we can use dominated convergence to conclude that $\lim_{s\downarrow 0} \mathbb{E}[\mathrm{e}^{sX'_1}] = 1$. Then, we have

$$\begin{aligned}
\lim_{s\downarrow 0} s^{-1} \log \mathbb{E}[\mathrm{e}^{sX'_1}] &= \lim_{s\downarrow 0} s^{-1} \log(1 + (\mathbb{E}[\mathrm{e}^{sX'_1}] - 1)) \\
&= \lim_{s\downarrow 0} s^{-1}(\mathbb{E}[\mathrm{e}^{sX'_1}] - 1) \\
&= \mathbb{E}[\lim_{s\downarrow 0} s^{-1}(\mathrm{e}^{sX'_1} - 1)] \\
&= \mathbb{E}[X'_1] \\
&= \kappa^{-1}.
\end{aligned}$$

Here, the third equality also follows from dominated convergence using the bound $s^{-1}(\mathrm{e}^{sX'_1} - 1) < s_0^{-1}(\mathrm{e}^{s_0X'_1} - 1)$ for all $0 < s \leq s_0$, where $\mathbb{E}[s_0^{-1}(\mathrm{e}^{s_0X'_1} - 1)] < \infty$. Therefore, there exists $c'_1 > 0$ such that

$$\limsup_{t\uparrow\infty} t^{-1} \log \mathbb{P}\left(\sum_{i=1}^{\lceil t(\kappa-\varepsilon)\rceil} X'_i \geq t\right) < -c'_1. \tag{16}$$

Again, by using the exponential Markov inequality on the other term, we obtain, for all $s > 0$,

$$\limsup_{t\uparrow\infty} t^{-1} \log \mathbb{P}\left(\sum_{i=1}^{\lfloor t(\kappa+\varepsilon)\rfloor} -X'_i > -t\right) \leq s(1+(\kappa+\varepsilon)s^{-1}\log \mathbb{E}[\mathrm{e}^{-sX'_1}]).$$

Since, by the dominated-convergence theorem, $\lim_{s\downarrow 0}\mathbb{E}[e^{-sX'_1}]=1$, we have

$$\begin{aligned}\lim_{s\downarrow 0} s^{-1}\log\mathbb{E}[\mathrm{e}^{-sX'_1}] &= \lim_{s\downarrow 0} s^{-1}\log\left(1-(1-\mathbb{E}[\mathrm{e}^{-sX'_1}])\right)\\ &= \lim_{s\downarrow 0} -s^{-1}(1-\mathbb{E}[\mathrm{e}^{-sX'_1}])\\ &= \mathbb{E}[\lim_{s\downarrow 0} -s^{-1}(1-\mathrm{e}^{-sX'_1})]\\ &= \mathbb{E}[-X'_1]\\ &= -\kappa^{-1}.\end{aligned}$$

Here, the third equality follows from the dominated-convergence theorem, using the bound $s^{-1}\big(1-\mathrm{e}^{-sX'_1}\big)<X'_1$, where $X'_1$ has finite expectation. Therefore, there exists $c'_2>0$ such that

$$\limsup_{t\uparrow\infty} t^{-1}\log\mathbb{P}\left(\sum_{i=1}^{\lfloor t(\kappa+\varepsilon)\rfloor} -X'_i > -t\right) < -c'_2. \tag{17}$$

Now, combining (15), (16), and (17) proves the lemma.

*Proof of Proposition* 3. Let us fix $\varepsilon'>0$ and note that

$$\begin{aligned}\mathbb{P}(|\mathcal{Y}^\kappa_t-\mathcal{Y}'_t|\ge\delta t^{1/2+\varepsilon}) \le{}& \mathbb{P}\left(\left\{\left|\sum_{i=\lfloor t\kappa\rfloor+1}^{K'_t} Y'_i\right|\ge\delta t^{1/2+\varepsilon}\right\}\cap\left\{\kappa\le\frac{K'_t}{t}\le\kappa+\varepsilon'\right\}\right)\\ &+\mathbb{P}\left(\left\{\left|\sum_{i=K'_t+1}^{\lfloor t\kappa\rfloor} Y'_i\right|\ge\delta t^{1/2+\varepsilon}\right\}\cap\left\{\kappa-\varepsilon'\le\frac{K'_t}{t}\le\kappa\right\}\right)\\ &+\mathbb{P}\left(\frac{K'_t}{t}\notin B_{\varepsilon'}(\kappa)\right),\end{aligned} \tag{18}$$

where, by Lemma 8, $\limsup_{t\uparrow\infty} t^{-2\varepsilon}\log\mathbb{P}(K'_t/t\notin B_{\varepsilon'}(\kappa))=-\infty$ and, hence, the third summand plays no role, by logarithmic equivalence. For the first summand in (18), since the $Y'_i$ are i.i.d., we have the bound

$$\sum_{m=\lfloor t\kappa\rfloor+1}^{\lfloor t(\kappa+\varepsilon')\rfloor}\mathbb{P}\left(\left\{\left|\sum_{i=\lfloor t\kappa\rfloor+1}^{K'_t} Y'_i\right|\ge\delta t^{1/2+\varepsilon}\right\}\cap\{K'_t=m\}\right)\le\sum_{m=1}^{r_1(t,\varepsilon')}\mathbb{P}\left(\left|\sum_{i=1}^{m} Y'_i\right|\ge\delta t^{1/2+\varepsilon}\right), \tag{19}$$

where the number of summands in the above sum is $r_1(t,\varepsilon'):=\lfloor t(\kappa+\varepsilon')\rfloor-\lfloor t\kappa\rfloor$. Note that the corresponding upper bound is valid also for the second summand in (18) with $r_1(t,\varepsilon')$ replaced by $r_2(t,\varepsilon'):=\lfloor t\kappa\rfloor-\lfloor t(\kappa-\varepsilon')\rfloor$. Furthermore, since $\limsup_{t\uparrow\infty} t^{-2\varepsilon}\log r(t,\varepsilon')=0$, where $r(t,\varepsilon')=r_1(t,\varepsilon')\vee r_2(t,\varepsilon')$, we have

$$\begin{aligned}\limsup_{t\uparrow\infty} t^{-2\varepsilon}\log\mathbb{P}(|\mathcal{Y}^\kappa_t-\mathcal{Y}'_t|\ge\delta t^{1/2+\varepsilon}) &\le \limsup_{t\uparrow\infty} t^{-2\varepsilon}\sup_{0\le\alpha\le\varepsilon'}\log\mathbb{P}\left(\left|\sum_{i=1}^{\lfloor t\alpha\rfloor} Y'_i\right|\ge\delta t^{1/2+\varepsilon}\right)\\ &\le \limsup_{t\uparrow\infty} t^{-2\varepsilon}\log\mathbb{P}\left(\left|\sum_{i=1}^{\lfloor t\varepsilon'\rfloor} Y'_i\right|\ge\delta t^{1/2+\varepsilon}\right),\end{aligned}$$

where we used the fact that, by symmetry and independence,

$$\mathbb{P}\left(\sum_{i=1}^{\lfloor t\alpha\rfloor} Y_i' \geq \delta t^{1/2+\varepsilon}\right) = 2\mathbb{P}\left(\sum_{i=1}^{\lfloor t\alpha\rfloor} Y_i' \geq \delta t^{1/2+\varepsilon}, \sum_{i=\lfloor t\alpha\rfloor+1}^{\lfloor t\varepsilon'\rfloor} Y_i' \geq 0\right) \leq 2\mathbb{P}\left(\sum_{i=1}^{\lfloor t\varepsilon'\rfloor} Y_i' \geq \delta t^{1/2+\varepsilon}\right),$$

and similar for $\mathbb{P}(\sum_{i=1}^{\lfloor t\alpha\rfloor} Y_i' \leq -\delta t^{1/2+\varepsilon})$. Hence, using the moderate-deviation principle, Proposition 2, we have

$$\limsup_{t\uparrow\infty} t^{-2\varepsilon} \log \mathbb{P}(|\mathcal{Y}_t^{\kappa} - \mathcal{Y}_t'| \geq \delta t^{1/2+\varepsilon}) \leq -I_{\lambda,\theta}^{\varepsilon'}(\delta),$$

which tends to $-\infty$ as $\varepsilon'$ tends to 0, as desired.

## 6. Proofs for the Moderate Deviations in Case $\theta = \pi/2$

As we can see, for $\theta = \pi/2$, the history set will never vanish, so we never observe independence between consecutive steps. In this case, we define the renewal steps as in [6, Section 4]. For a constant $\varkappa \geq 6$ and $u \in \mathbb{R}^2$, we define $u^{\rightarrow} := u + (\varkappa, 0)$. We set $\varsigma_0 = 0$ and then iteratively define

$$\varsigma_j := \inf \left\{ n > \varsigma_{j-1} : \begin{array}{c} L_n \leq \varkappa, \sum_{i=\varsigma_{j-1}+1}^{n} X_i \geq \varkappa + 1, \\ \#\left(\mathcal{C}_{\pi/2}(V_n) \cap B(V_n, \varkappa+1) \cap \mathcal{P}_\lambda\right) = 1 \text{ and} \\ \#\left(\mathcal{C}_{\pi/2}(V_n^{\rightarrow}) \cap B(V_n^{\rightarrow}, 1) \cap \mathcal{P}_\lambda\right) = 1 \end{array} \right\}.$$

We call the $\varsigma_j$th step the $j$th renewal step. At a renewal step, the width of the history set is at most $\varkappa$, and the condition $\sum_{i=\varsigma_{j-1}+1}^{\varsigma_j} X_i \geq \varkappa + 1$ ensures that the history sets at different renewal steps are disjoint. There is only one Poisson point in the set $\mathcal{C}_{\pi/2}(V_{\varsigma_j}) \cap B(V_{\varsigma_j}, \varkappa+1)$, which is actually included in the subset $\mathcal{C}_{\pi/2}(V_{\varsigma_j}^{\rightarrow}) \cap B(V_{\varsigma_j}^{\rightarrow}, 1)$. Therefore, this Poisson point is the next point in our exploration, i.e. it is $V_{\varsigma_j+1}$. Moreover, if we had started our exploration from $V_{\varsigma_j}^{\rightarrow}$ then $V_{\varsigma_j+1}$ would still have been the next exploration point. As noted in [6], conditional on being at the $j$th renewal step, this Poisson point $V_{\varsigma_j+1}$ is uniformly distributed on $\mathcal{C}_{\pi/2}(V_{\varsigma_j}^{\rightarrow}) \cap B(V_{\varsigma_j}^{\rightarrow}, 1)$. Thus, for $j \geq 1$, the paths starting from $V_{\varsigma_j}^{\rightarrow}$ and ending at $V_{\varsigma_{j+1}}$ are i.i.d. copies of the path starting at $o$ until the first renewal step, with the following initial conditions.

(a) A single point is uniformly distributed in $\mathcal{C}_{\pi/2}(o) \cap B(o, 1)$.

(b) The set $\left(\mathcal{C}_{\pi/2}((-\varkappa, 0)) \cap B((-\varkappa, 0), \varkappa+1)\right) \setminus \left(\mathcal{C}_{\pi/2}(o) \cap B(o, 1)\right)$ contains no points.

(c) An independent Poisson point process is placed on $\left(\mathcal{C}_{\pi/2}((-\varkappa, 0)) \cap B((-\varkappa, 0), \varkappa+1)\right)^c$.

Let $\mathfrak{Z} = (\mathfrak{Z}_1, \mathfrak{Z}_2)$ be the position of the path at the first renewal step. Then, writing

$$\mathfrak{U}_{j+1}' := V_{\varsigma_{j+1}} - V_{\varsigma_j},$$

the sequence $\{\mathfrak{U}_{j+1}'\}_{j\geq 1}$ consists of i.i.d. copies of $\mathfrak{Z} + (\varkappa, 0)$, as stated in [6, Proposition 4.5]. Furthermore, as in [6, Proposition 4.2], there exist constants $\mathfrak{c}, \mathfrak{C} > 0$ such that, for any $j \geq 0$ and any $n \geq 1$,

$$\mathbb{P}(\varsigma_{j+1} - \varsigma_j > n) \leq \mathfrak{C}\mathrm{e}^{-\mathfrak{c}n}. \tag{20}$$

As earlier, we now define

$$\mathfrak{K}'_t := \sup\{n > 0 \colon \varsigma_n \le K_t\} \tag{21}$$

as the index of the largest stopping time before $K_t$. In particular, writing

$$\mathcal{Y}''_t := \sum_{i=2}^{\mathfrak{K}'_t} \mathfrak{Y}'_i,$$

where $\mathfrak{U}'_i =: (\mathfrak{X}'_i, \mathfrak{Y}'_i)$ in Cartesian coordinates, we have, for $\mathfrak{K}'_t \ge 1$,

$$\sum_{i=1}^{K_t} Y_i = \sum_{i=1}^{\mathfrak{K}'_t} \mathfrak{Y}'_i + \sum_{i=\varsigma_{\mathfrak{K}'_t}+1}^{K_t} Y_i = \mathfrak{Y}'_1 + \mathcal{Y}''_t + \sum_{i=\varsigma_{\mathfrak{K}'_t}+1}^{K_t} Y_i. \tag{22}$$

Note that, since $\mathfrak{Y}'_1$ does not have the same distribution as the other $\mathfrak{Y}'_i$s, we have taken the sum from $i = 2$ in the definition of $\mathcal{Y}''_t$ to ensure that $\mathcal{Y}''_t$ is a sum of i.i.d. random variables. As earlier, for the moderate deviations, the difference between $\mathcal{Y}_t$ and $\mathcal{Y}''_t$ is irrelevant in the following sense.

**Proposition 5.** *For any $\lambda > 0$ and $0 < \varepsilon < 1/2$, $\{t^{-1/2-\varepsilon}\mathcal{Y}_t\}_{t\ge 0}$ and $\{t^{-1/2-\varepsilon}\mathcal{Y}''_t\}_{t\ge 0}$ are exponentially equivalent, i.e. for any $\delta > 0$,*

$$\limsup_{t\uparrow\infty} t^{-2\varepsilon} \log \mathbb{P}(|\mathcal{Y}_t - \mathcal{Y}''_t| \ge \delta t^{1/2+\varepsilon}) = -\infty.$$

Note that we cannot use Lemma 5 to obtain a bound for the $Y_i$s in the case $\theta = \pi/2$. Therefore, we need the following result to prove Proposition 5.

**Lemma 9.** *For any $0 < \theta \le \pi/2$, there exist two sequences of nonnegative random variables $\{\tilde{X}_i\}_{i\ge 1}$ and $\{\tilde{Y}_i\}_{i\ge 1}$, defined on the same probability space as $\{(X_i, Y_i)\}_{i\ge 1}$, such that*

(i) *$\{\tilde{X}_i\}_{i\ge 1}$ are i.i.d. with $\mathbb{P}(\tilde{X}_1 \ge t) = \mathrm{e}^{-\lambda\theta t^2}$ for all $t \ge 0$, and $\sum_{i=1}^n X_i \le \sum_{i=1}^n \tilde{X}_i$ for all $n \ge 1$ a.s.;*

(ii) *$\{\tilde{Y}_i\}_{i\ge 1}$ are i.i.d. with $\mathbb{P}(\tilde{Y}_1 \ge t) = \mathrm{e}^{-\lambda\theta t^2/2}$ for all $t \ge 0$, and $\sum_{i=1}^n Y_i \le \sum_{i=1}^n \tilde{Y}_i$ for all $n \ge 1$ a.s.*

*Proof of Lemma 9.* We treat the domination in the vertical and horizontal direction separately. Let us start with the horizontal direction.

*Part (i): Existence of $\{\tilde{X}_i\}_{i\ge 0}$.*Let $\{U_i\}_{i\ge 1}$ be the sequence of progress random variables. Recall that we write $U_i = (X_i, Y_i)$ in Cartesian coordinates and $U_i = (R_i, \Phi_i)$ in polar coordinates. Recall further that $V_0 = o$ and $V_k = (\sum_{i=1}^k X_i, \sum_{i=1}^k Y_i)$ for all $k \ge 1$. Now, we set $\tilde{V}_0 := V_0$ and define the sequences $\{\tilde{V}_k\}_{k\ge 1}$, $\{\tilde{W}_k\}_{k\ge 1}$, and $\{\tilde{X}_k\}_{k\ge 1}$ using the following recursive equations:

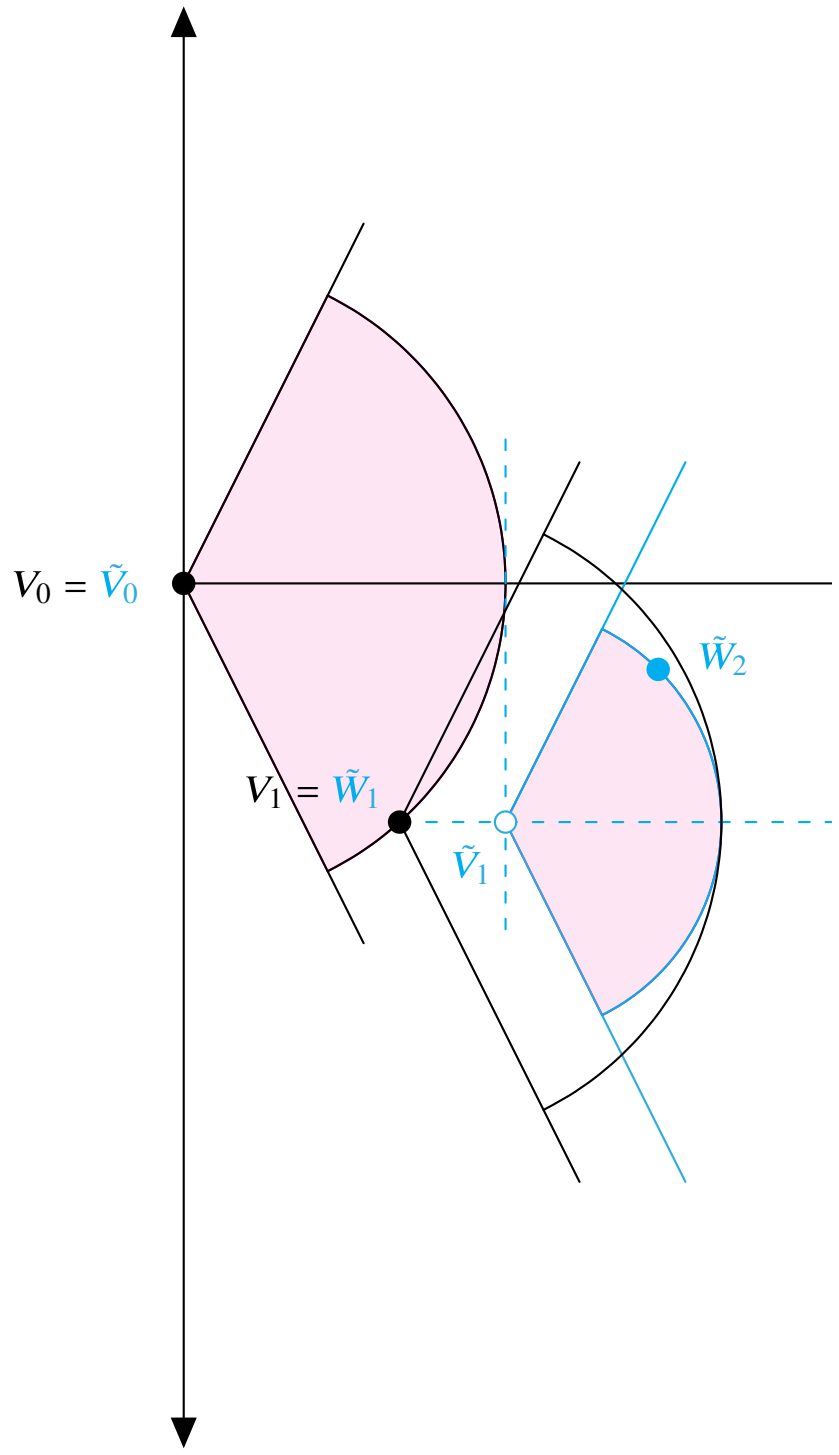

FIGURE 3. A realization of part (i): existence of $\{\tilde{X}_i\}_{i\geq 0}$.

$$\tilde{W}_k := \arg\min \left\{ \left| v - \tilde{V}_{k-1} \right| : v \in \mathcal{P} \cap C_\theta(\tilde{V}_{k-1}) \right\},$$
$$\tilde{X}_k := |\tilde{W}_k - \tilde{V}_{k-1}| \text{ and}$$
$$\tilde{V}_k := \left( \sum_{i=1}^{k} \tilde{X}_i, \sum_{i=1}^{k} Y_i \right)$$

(see Figure 3 for an illustration). Namely, the sequence of new points $\{\tilde{V}_k\}_{k\geq 1}$ is constructed in such a way that it dominates the original navigation along the $x$ axis, by considering $\tilde{X}_k$ to be the maximum progress made along the $x$ axis and the first coordinate of $\tilde{V}_k$ to be the total maximum progress made by the original navigation along the $x$ axis.

Note that, by construction, for all $k \geq 1$ and for all $1 \leq i \leq k-1$,

$$C_\theta(\tilde{V}_k) \cap B(\tilde{V}_i, \tilde{X}_{i+1}) = \emptyset.$$

As a consequence, by the total independence of the underlying Poisson point process, $\{\tilde{X}_k\}_{k\geq 1}$ are i.i.d. Also, by construction,

$$\mathbb{P}(\tilde{X}_1 \geq t) = \mathrm{e}^{-\lambda\theta t^2}.$$

Now, we show the domination by induction. First, note that $X_1 \leq R_1 = \tilde{X}_1$. Then, by assuming that $\sum_{i=1}^{k} X_i \leq \sum_{i=1}^{k} \tilde{X}_i$, we note that

$$C_\theta(V_k) \supseteq C_\theta(\tilde{V}_k) \ni \tilde{W}_{k+1}.$$

Therefore,

$$X_{k+1} \le R_{k+1} \le |\tilde{W}_{k+1} - V_k| \le |\tilde{W}_{k+1} - \tilde{V}_k| + |\tilde{V}_k - V_k| = \tilde{X}_{k+1} + \sum_{i=1}^{k} \tilde{X}_i - \sum_{i=1}^{k} X_i,$$

which implies that $\sum_{i=1}^{k+1} X_i \le \sum_{i=1}^{k+1} \tilde{X}_i$. This, by induction, proves Lemma 9(i).

*Part (ii): Existence of* $\{\tilde{Y}_i\}_{i\ge 0}$. Given an angle $\theta > 0$, we define

$$C_\theta^u := \{y = (r_y, \varphi_y) : r_y > 0, 0 \le \varphi_y \le \theta\}$$

and, for any $v \in \mathbb{R}^2$, we write $C_\theta^u(v) := v + C_\theta^u$. Then, similar to part (i), we set $\hat{V}_0 := V_0$ and define sequences $\{\hat{V}_k\}_{k\ge1}$, $\{\hat{W}_k\}_{k\ge1}$, and $\{\tilde{Y}_k\}_{k\ge1}$ using the recursive equations

$$\hat{W}_k := \arg\min\left\{\left|v - \hat{V}_{k-1}\right| : v \in \mathcal{P} \cap C_\theta^u(\hat{V}_{k-1})\right\},$$

$$\tilde{Y}_k := |\hat{W}_k - \hat{V}_{k-1}|, \text{ and}$$

$$\hat{V}_k := \left(\sum_{i=1}^{k} X_i - \sum_{i=1}^{k} Y_i \cot\theta + \sum_{i=1}^{k} \tilde{Y}_i \cos\theta, \sum_{i=1}^{k} \tilde{Y}_i \sin\theta\right);$$

see Figure 4 for an illustration. Similar to above, $\{\sum_{i=1}^{k} \tilde{Y}_i \sin\theta\}_{k\ge1}$, the sequence of second coordinates of $\hat{V}_k$, makes sure that the navigation stays dominated along the $y$ axis, as $\tilde{Y}_k \sin\theta$ is the maximum progress that can be made along the $y$ axis starting from $\hat{V}_{k-1}$. Whereas the first coordinates of $\hat{V}_k$ are chosen in such a way that the domination holds, as described in (23). Note that, by construction, for all $k \ge 1$ and for all $1 \le i \le k-1$,

$$C_\theta^u(\hat{V}_k) \cap B(\hat{V}_i, \tilde{Y}_{i+1}) = \emptyset.$$

Hence, by the independence of the Poisson point process in disjoint regions, $\{\tilde{Y}_k\}_{k\ge1}$ are i.i.d. Also, by construction,

$$\mathbb{P}(\tilde{Y}_1 \ge t) = \mathrm{e}^{-\lambda\theta t^2/2}.$$

Now, we show the domination by induction. First, note that $Y_1 \le R_1 \sin\theta \le \tilde{Y}_1 \sin\theta$. Then, by assuming that $\sum_{i=1}^{k} Y_i \le \sum_{i=1}^{k} \tilde{Y}_i \sin\theta$, we note that

$$C_\theta(V_k) \supseteq C_\theta^u(\hat{V}_k) \ni \hat{W}_{k+1}.$$

Therefore,

$$\begin{aligned}\frac{Y_{k+1}}{\sin\theta} &\le R_{k+1} \le |\hat{W}_{k+1} - V_k| \le |\hat{W}_{k+1} - \hat{V}_k| + |\hat{V}_k - V_k| \\ &= \tilde{Y}_{k+1} + \frac{1}{\sin\theta}\left(\sum_{i=1}^{k} \tilde{Y}_i \sin\theta - \sum_{i=1}^{k} Y_i\right), \end{aligned} \tag{23}$$

which implies that $\sum_{i=1}^{k+1} Y_i \le \sum_{i=1}^{k+1} \tilde{Y}_i \sin\theta$. This, by induction, proves Lemma 9(ii), thus completing the proof.

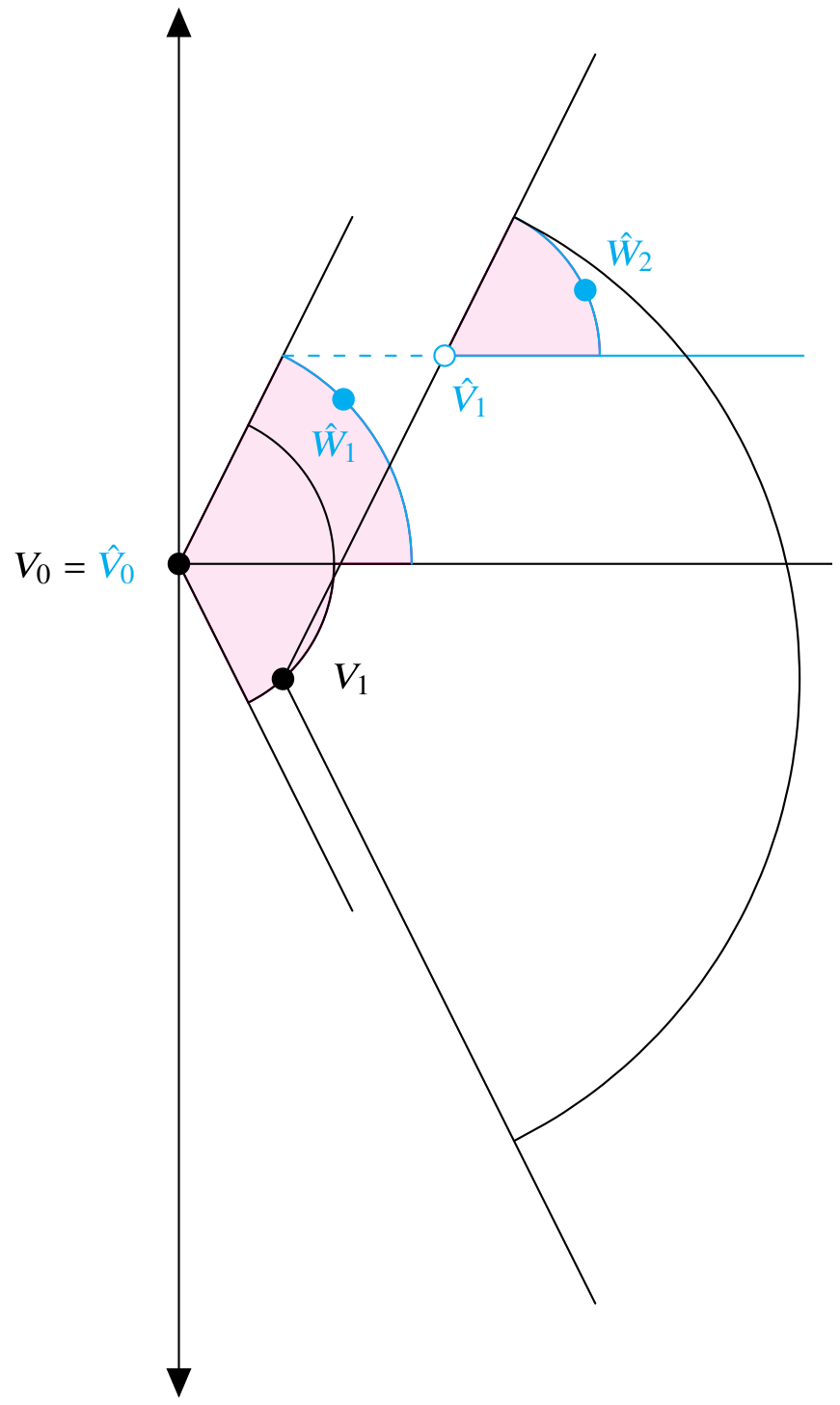

FIGURE 4. A realization of part (ii): existence of $\{\tilde{Y}_i\}_{i\geq 0}$.

*Proof of Proposition* 5. We can proceed as in the proof of Proposition 1. Note that, by (22),

$$|\mathcal{Y}_t - \mathcal{Y}_t''| \leq \left|\sum_{i=1}^{\varsigma_1} Y_i\right| + \max\left\{\left|\sum_{i=\varsigma_{\mathfrak{K}_t'}+1}^{K_t} Y_i\right|, \left|\sum_{i=\varsigma_{\mathfrak{K}_t'}+1}^{K_t+1} Y_i\right|\right\}$$

and, therefore, by symmetry,

$$\begin{aligned}
\mathbb{P}(|\mathcal{Y}_t - \mathcal{Y}_t''| \geq \delta t^{1/2+\varepsilon}) &\leq \mathbb{P}\left(\left|\sum_{i=1}^{\varsigma_1} Y_i\right| \geq \frac{1}{2}\delta t^{1/2+\varepsilon}\right) + \mathbb{P}\left(\left|\sum_{i=\varsigma_{\mathfrak{K}_t'}+1}^{K_t} Y_i\right| \geq \frac{1}{2}\delta t^{1/2+\varepsilon}\right) \\
&\quad + \mathbb{P}\left(\left|\sum_{i=\varsigma_{\mathfrak{K}_t'}+1}^{K_t+1} Y_i\right| \geq \frac{1}{2}\delta t^{1/2+\varepsilon}\right) \\
&= 2\mathbb{P}\left(\sum_{i=1}^{\varsigma_1} Y_i \geq \frac{1}{2}\delta t^{1/2+\varepsilon}\right) + 2\mathbb{P}\left(\sum_{i=\varsigma_{\mathfrak{K}_t'}+1}^{K_t} Y_i \geq \frac{1}{2}\delta t^{1/2+\varepsilon}\right) \\
&\quad + 2\mathbb{P}\left(\sum_{i=\varsigma_{\mathfrak{K}_t'}+1}^{K_t+1} Y_i \geq \frac{1}{2}\delta t^{1/2+\varepsilon}\right),
\end{aligned}$$

which, together with Lemma 9, yields

$$
\begin{aligned}
&\mathbb{P}(|\mathcal{Y}_t - \mathcal{Y}_t''| \geq \delta t^{1/2+\varepsilon}) \\
&\quad \leq 2\mathbb{P}\Bigg(\sum_{i=1}^{\varsigma_1} \tilde{Y}_i \geq \frac{1}{2}\delta t^{1/2+\varepsilon}\Bigg) + 4\mathbb{P}\Bigg(\sum_{i=\varsigma_{\mathfrak{K}_t'}+1}^{\varsigma_{\mathfrak{K}_t'+1}} \tilde{Y}_i \geq \frac{1}{2}\delta t^{1/2+\varepsilon}\Bigg) \\
&\quad \leq 2\mathbb{P}\Bigg(\sum_{i=1}^{\varsigma_1} \tilde{Y}_i \geq \frac{1}{2}\delta t^{1/2+\varepsilon}\Bigg) + 4\sum_{j=0}^{\lfloor t^2 \rfloor} \mathbb{P}\Bigg(\sum_{i=\varsigma_j+1}^{\varsigma_{j+1}} \tilde{Y}_i \geq \frac{1}{2}\delta t^{1/2+\varepsilon}, \mathfrak{K}_t' = j\Bigg) + 4\mathbb{P}(\mathfrak{K}_t' > \lfloor t^2 \rfloor) \\
&\quad \leq 6\mathbb{P}\Bigg(\sum_{i=1}^{\varsigma_1} \tilde{Y}_i \geq \frac{1}{2}\delta t^{1/2+\varepsilon}\Bigg) + 4\lfloor t^2 \rfloor \mathbb{P}\Bigg(\sum_{i=\varsigma_1+1}^{\varsigma_2} \tilde{Y}_i \geq \frac{1}{2}\delta t^{1/2+\varepsilon}\Bigg) + 4\mathbb{P}(\mathfrak{K}_t' > \lfloor t^2 \rfloor) \\
&\quad \leq (6 + 4\lfloor t^2 \rfloor)\mathbb{P}\Bigg(\sum_{i=1}^{\lfloor t^{\varepsilon'} \rfloor} \tilde{Y}_i \geq \frac{1}{2}\delta t^{1/2+\varepsilon}\Bigg) + 6\mathbb{P}(\varsigma_1 > \lfloor t^{\varepsilon'} \rfloor) + 4\lfloor t^2 \rfloor \mathbb{P}(\varsigma_2 - \varsigma_1 > \lfloor t^{\varepsilon'} \rfloor) \\
&\qquad + 4\mathbb{P}\Bigg(\sum_{i=1}^{\lfloor t^2 \rfloor} \mathfrak{X}_i' < t\Bigg)
\end{aligned}
$$

for some $\varepsilon' \in (2\varepsilon, \frac{1}{2} + \varepsilon)$. This, together with (20) and the exponential Markov inequality, implies that

$$
\begin{aligned}
\mathbb{P}(|\mathcal{Y}_t - \mathcal{Y}_t''| \geq \delta t^{1/2+\varepsilon}) &\leq (6 + 4\lfloor t^2 \rfloor)\mathrm{e}^{-\frac{1}{2}\delta t^{1/2+\varepsilon}} \mathbb{E}[\mathrm{e}^{\tilde{Y}_1}]^{\lfloor t^{\varepsilon'} \rfloor} \\
&\quad + (6 + 4\lfloor t^2 \rfloor)\mathfrak{C}\mathrm{e}^{-\mathfrak{c}\lfloor t^{\varepsilon'} \rfloor} + 4\mathrm{e}^t \mathbb{E}[\mathrm{e}^{-\mathfrak{X}_1'}]^{\lfloor t^2 \rfloor},
\end{aligned}
$$

and therefore,

$$
\begin{aligned}
&\limsup_{t\uparrow\infty} t^{-2\varepsilon} \log \mathbb{P}(|\mathcal{Y}_t - \mathcal{Y}_t''| \geq \delta t^{1/2+\varepsilon}) \\
&\quad \leq -\min\Big\{\liminf_{t\uparrow\infty} t^{-2\varepsilon}\big(\tfrac{1}{2}\delta t^{1/2+\varepsilon} - \lfloor t^{\varepsilon'} \rfloor \log \mathbb{E}[\mathrm{e}^{\tilde{Y}_1}]\big), \liminf_{t\uparrow\infty} c\lfloor t^{\varepsilon'} \rfloor t^{-2\varepsilon}, \\
&\qquad\qquad -\limsup_{t\uparrow\infty} t^{-2\varepsilon}(t + \lfloor t^2 \rfloor \log \mathbb{E}[\mathrm{e}^{-\mathfrak{X}_1'}])\Big\} \\
&\quad = -\infty,
\end{aligned}
$$

as desired.

*Proof of Theorem 1, case* $\theta = \pi/2$. Using Proposition 5, it suffices to consider $\mathcal{Y}_t''$. However, an argument similar to that in the proof of Propositions 2 and 3 gives us that, for any $0 < \lambda$ and $0 < \varepsilon < 1/2$, the process $\{t^{-1/2-\varepsilon}\mathcal{Y}_t''\}_{t\geq 0}$ obeys the moderate-deviation principle with rate function

$$
I_{\lambda,\pi/2}(x) := x^2 \frac{\mathbb{E}[\mathfrak{Z}_1]}{2\mathbb{E}[\mathfrak{Z}_2^2]},
$$

and hence, setting

$$\rho(\lambda, \pi/2) := \frac{\mathbb{E}[\mathfrak{Z}_1]}{2\mathbb{E}[\mathfrak{Z}_2^2]}$$

gives us the required form. From the definition of moderate deviations, it follows that, by scaling both the coordinates by $\sqrt{\lambda}$, we have

$$I_{\lambda,\pi/2}(x) = \lambda^{\varepsilon} I_{1,\pi/2}(x\lambda^{1/4-\varepsilon/2}),$$

and hence,

$$\rho(\lambda, \pi/2) = \sqrt{\lambda}\rho(1, \pi/2).$$

This completes the proof of the theorem.

## 7. Proofs for the Large Deviations

### 7.1. The independent case $0 < \theta \le \pi/4$

*Proof of Lemma* 1. The statement follows from an application of the multivariate Cramér theorem for empirical means of sequences of i.i.d. random variables that possess exponential moments; see [7, Corollary 6.1.6]. In order to compute the exponential moments, we use polar coordinates, i.e. consider $U_1 = (R, \Phi) \in \mathbb{R}_+ \times [-\pi, \pi)$. As already used earlier, the radius follows a Rayleigh distribution, i.e.

$$\mathbb{P}(R > r) = \exp(-\lambda\theta r^2),$$

and we note that, due to the isotropy of the model, $\Phi$ is uniformly distributed in $[-\theta, \theta]$. Hence,

$$\begin{aligned}\mathbb{E}[\exp(\langle\gamma, U_1\rangle)] &= \int_0^\infty \mathrm{d}r \frac{1}{2\theta}\int_{-\theta}^{\theta} \mathrm{d}\varphi\ \exp(\gamma_1 r\cos\varphi + \gamma_2 r\sin\varphi - \lambda\theta r^2)\lambda 2\theta r \\ &= \lambda\int_0^\infty \mathrm{d}r\, r\exp(-\lambda\theta r^2)\int_{-\theta}^{\theta}\mathrm{d}\varphi\ \exp(\gamma_1 r\cos\varphi + \gamma_2 r\sin\varphi),\end{aligned} \tag{24}$$

as desired. In particular,

$$\mathrm{e}^{J_{\lambda,\theta}(\gamma)} \le 2\pi\lambda\int_0^\infty \mathrm{d}r\, r\exp(-\lambda\theta r^2)\exp(r(|\gamma_1| + |\gamma_2|)) < \infty.$$

Since $J_{\lambda,\theta}$ is strictly convex and differentiable, by [14, Theorem 1], its Legendre transform $\mathcal{J}_{\lambda,\theta}$ is strictly convex and differentiable on $\{u \in \mathbb{R}^2 : \mathcal{J}_{\lambda,\theta}(u) < \infty\}$.

Further note that, writing $u = (r_u, \varphi_u)$ in polar coordinates,

$$\mathrm{e}^{-\mathcal{J}_{\lambda,\theta}(u)} = \inf_{q\ge 0,\, \alpha\in[-\pi,\pi]}\int_0^\infty \mathrm{d}r\int_{-\theta}^{\theta}\mathrm{d}\varphi\, \lambda r\exp(-\lambda\theta r^2 + q[r\cos(\varphi-\alpha) - r_u\cos(\varphi_u-\alpha)]),$$

and, if $\varphi_u \notin (-\theta, \theta)$, then we can always choose $\alpha_0 \in [-\pi, \pi]$ such that, for all $\varphi \in [-\theta, \theta]$, we have $r\cos(\varphi - \alpha_0) - r_u\cos(\varphi_u - \alpha_0) \le 0$. One choice of $\alpha_0$ for example is given by

$$\alpha_0 = \begin{cases}\pi/2 + \theta & \text{if } \varphi_u \in [\theta, \pi], \\ -\pi/2 - \theta & \text{if } \varphi_u \in [-\pi, -\theta].\end{cases}$$

Then, since

$$\int_0^\infty \mathrm{d}r \int_{-\theta}^{\theta} \mathrm{d}\varphi\, \lambda r \exp(-\lambda\theta r^2) < \infty,$$

by the dominated-convergence theorem, we obtain

$$\begin{aligned} \mathrm{e}^{-\mathcal{J}_{\lambda,\theta}(u)} &\le \lim_{q\to\infty} \int_0^\infty \mathrm{d}r \int_{-\theta}^{\theta} \mathrm{d}\varphi\, \lambda r \exp(-\lambda\theta r^2 + q[r\cos(\varphi-\alpha_0) - r_u \cos(\varphi_u - \alpha_0)]) \\ &= \int_0^\infty \mathrm{d}r \int_{-\theta}^{\theta} \mathrm{d}\varphi \lim_{q\to\infty} \lambda r \exp(-\lambda\theta r^2 + q[r\cos(\varphi-\alpha_0) - r_u \cos(\varphi_u - \alpha_0)]) \\ &= 0. \end{aligned}$$

Therefore, we get $\mathcal{J}_{\lambda,\theta}(u) = \infty$ whenever $\varphi_u \notin (-\theta, \theta)$. A similar calculation also gives us that $\mathcal{J}_{\lambda,\theta}(u) = \infty$ for $r_u = 0$. This implies that $\mathcal{J}_{\lambda,\theta}(u) = \infty$ for $u \notin \mathcal{C}^o_\theta$.

Now, we show that $\mathcal{J}_{\lambda,\theta}(u) < \infty$ for $u \in \mathcal{C}^o_\theta$. In this case, $r_u > 0$ and $\varphi_u \in (-\theta, \theta)$. For this, we pick $0 < \eta < \pi/16$ small enough such that $[\varphi_u - 4\eta, \varphi_u + 4\eta] \in (-\theta, \theta)$ and define two sets $A_\alpha$, $B_\alpha$ as

$$A_\alpha := \begin{cases} (r_u/(\sin\eta), \infty) & \text{if } |\alpha - \varphi_u| \le \pi/2 + \eta, \\ (0, r_u \sin\eta) & \text{if } |\alpha - \varphi_u| > \pi/2 + \eta, \end{cases}$$

and

$$B_\alpha := \begin{cases} [\varphi_u + 2\eta, \varphi_u + 3\eta] & \text{if } \alpha \in [\varphi_u, \varphi_u + \pi], \\ [\varphi_u - 3\eta, \varphi_u - 2\eta] & \text{if } \alpha \in [\varphi_u - \pi, \varphi_u). \end{cases}$$

Here, we consider the range of $\alpha$ as $[\varphi_u - \pi, \varphi_u + \pi]$ instead of $[-\pi, \pi]$ for simplicity. Now, for any $\alpha \in [\varphi_u - \pi, \varphi_u + \pi]$, $r \in A_\alpha$, and $\varphi \in B_\alpha$, we have $r\cos(\varphi - \alpha) - r_u \cos(\varphi_u - \alpha) \ge 0$. Therefore,

$$\begin{aligned} \mathrm{e}^{-\mathcal{J}_{\lambda,\theta}(u)} &\ge \inf_{\alpha\in[\varphi_u-\pi,\varphi_u+\pi]} \int_{A_\alpha} \mathrm{d}r \int_{B_\alpha} \mathrm{d}\varphi\, \lambda r \exp(-\lambda\theta r^2) \\ &= \eta\lambda \min\left\{ \int_0^{r_u \sin\eta} \mathrm{d}r\, r\exp(-\lambda\theta r^2), \int_{r_u/(\sin\eta)}^\infty \mathrm{d}r\, r\exp(-\lambda\theta r^2) \right\} \\ &> 0, \end{aligned}$$

which implies that, for all $u \in \mathcal{C}^o_\theta$, we have $\mathcal{J}_{\lambda,\theta}(u) < \infty$. This completes the proof.

In order to establish certain bounds, it will be convenient to also have a large-deviation result for i.i.d. sequences of progress variables $\overline{U}_i = (X_i, |Y_i|)$ at our disposal. For this, let $\overline{J}_{\lambda,\theta}$ denote the logarithmic moment-generating function of $\overline{U}_1$.

**Lemma 10.** *The sequence of i.i.d. copies* $\{\overline{U}_i\}_{i\ge 1}$ *of the progress variable* $\overline{U}_1 \in \mathbb{R} \times \mathbb{R}_+$ *satisfies the large-deviation principle with rate n and rate function*

$$\overline{\mathcal{J}}_{\lambda,\theta}(u) := \sup\{\langle \gamma, u\rangle - \overline{J}_{\lambda,\theta}(\gamma) \colon \gamma \in \mathbb{R}^2\},$$

*where* $\overline{J}_{\lambda,\theta}(\gamma) < \infty$ *for all* $\gamma \in \mathbb{R}^2$.

*Proof of Lemma* 10. Using the same arguments as in the proof of Lemma 1, we have

$$\mathbb{E}[\exp(\langle\gamma, \overline{U}_1\rangle)] = 2\lambda \int_0^\infty \mathrm{d}r\, r \exp(-\lambda\theta r^2) \int_0^\theta \mathrm{d}\varphi\, \exp(\gamma_1 r\cos\varphi + \gamma_2 r \sin\varphi) < \infty,$$

which, together with Cramér's theorem, gives the result.

*Proof of Theorem* 2. The key observation is that the sequence of progress variables $\{U_i\}_{i\geq 1}$ is i.i.d. for $\theta < \pi/4$. Indeed, note that, for all $i \geq 1$, $(C_\theta(V_i) \cap B(V_i, |U_{i+1}|)) \cap C_\theta(V_{i+1}) = \emptyset$. Hence, in every step, the navigation discovers a previous undiscovered region in space and the corresponding Poisson point clouds are stochastically independent. However, even though the progress steps are i.i.d., the statement is not a direct application of, e.g. Cramér's theorem, since $\{\mathcal{Y}_t\}_{t\geq 0}$ only tracks the vertical displacement; however, there is also a random horizontal displacement. Let us write $U_i = (X_i, Y_i)$, where $X_i$ is the first and $Y_i$ the second Cartesian coordinate of $U_i$ and recall that $K_t := \sup\{n > 0\colon \sum_{i=1}^n X_i < t\}$ denotes the number of steps the navigation takes before reaching $t$ along the $x$ axis.

***Step 1:*** *Lower bound for upper tail.* Let us start with the lower bound for the upper tail. It suffices to consider $a > 0$. Then, for all $1 > b > 0$, $c > 0$, and $\alpha > \beta > 0$,

$$\begin{aligned}
\mathbb{P}(\mathcal{Y}_t > at) &\geq \mathbb{P}(\mathcal{Y}_t > at,\, \lfloor\beta t\rfloor \leq K_t < \lfloor\alpha t\rfloor)\\
&\geq \mathbb{P}\Bigg(\sum_{i=1}^{\lfloor\beta t\rfloor} Y_i > at + \sum_{j=\lfloor\beta t\rfloor+1}^{\lfloor\alpha t\rfloor} |Y_j|,\, bt < \sum_{i=1}^{\lfloor\beta t\rfloor} X_i < t,\, \sum_{j=\lfloor\beta t\rfloor+1}^{\lfloor\alpha t\rfloor} X_j > (1-b)t\Bigg)\\
&\geq \mathbb{P}\Bigg(\sum_{j=\lfloor\beta t\rfloor+1}^{\lfloor\alpha t\rfloor} |Y_j| < ct,\, \sum_{j=\lfloor\beta t\rfloor+1}^{\lfloor\alpha t\rfloor} X_j > (1-b)t\Bigg)\mathbb{P}\Bigg(\sum_{i=1}^{\lfloor\beta t\rfloor} Y_i > (a+c)t,\, bt < \sum_{i=1}^{\lfloor\beta t\rfloor} X_i < t\Bigg)
\end{aligned}$$

by independence. Consequently, for $\delta = \alpha - \beta$,

$$\begin{aligned}
\liminf_{t\uparrow\infty} t^{-1} \log \mathbb{P}(\mathcal{Y}_t > at) \geq{}& -\delta \inf\{\overline{\mathcal{J}}_{\theta,\lambda}(x, y)\colon x > (1-b)/\delta,\, y < c/\delta\}\\
&- \beta \inf\{\mathcal{J}_{\theta,\lambda}(x, y)\colon b/\beta < x < 1/\beta,\, y > (a+c)/\beta\}.
\end{aligned}$$

Sending $b$ to 1 and fixing $\delta = \delta(c)$ such that $c/\delta > \mathbb{E}[|Y_1|]$, the first summand vanishes and we have

$$\liminf_{t\uparrow\infty} t^{-1} \log \mathbb{P}(\mathcal{Y}_t > at) \geq -\beta \mathcal{J}_{\theta,\lambda}(1/\beta, (a+c)/\beta).$$

Sending $c$ to 0, we arrive at

$$\liminf_{t\uparrow\infty} t^{-1} \log \mathbb{P}(\mathcal{Y}_t > at) \geq -\inf\{\beta \mathcal{J}_{\theta,\lambda}(1/\beta, a/\beta)\colon \beta > 0\}.$$

This expression reflects that the process has to find the optimal compromise between making the right number of steps for the displacement along the $y$ axis to be not too unlikely as well as hitting the time $t$ along the $x$ axis.

***Step 2:*** *Upper bound for upper tail.* For the upper bound, we can proceed similarly. For all $\alpha > 0$, we can bound

$$\begin{aligned}\mathbb{P}(\mathcal{Y}_t \geq at) &\leq \mathbb{P}\Bigg(\sum_{i=1}^{K_t} Y_i \geq at - |Y_{K_t+1}|, K_t \leq \lfloor \alpha t \rfloor\Bigg) + \mathbb{P}(K_t > \lfloor \alpha t \rfloor)\\ &= \sum_{m=0}^{\lfloor \alpha t \rfloor} \mathbb{P}\Bigg(\sum_{i=1}^{m} Y_i \geq at - |Y_{m+1}|, K_t = m\Bigg) + \mathbb{P}(K_t > \lfloor \alpha t \rfloor)\\ &\leq \sum_{m=0}^{\lfloor \alpha t \rfloor} \mathbb{P}\Bigg(\sum_{i=1}^{m} Y_i \geq t(a-\varepsilon), t(1-\varepsilon) \leq \sum_{i=1}^{m} X_i < t\Bigg) + \mathbb{P}(\mathfrak{B}_\varepsilon^c(t)) + \mathbb{P}(K_t > \lfloor \alpha t \rfloor),\end{aligned}$$

where $\mathfrak{B}_\varepsilon(t) = \{X_{m+1} \leq \varepsilon t, |Y_{m+1}| \leq \varepsilon t\}$ with $\varepsilon > 0$. Note that

$$\limsup_{t\uparrow\infty} t^{-1} \log \mathbb{P}(\mathfrak{B}_\varepsilon^c(t)) = -\infty,$$

since the upper tails of $|U_i|$ are of order $O(-t^2)$ on the exponential level. Furthermore,

$$\limsup_{\alpha\uparrow\infty} \limsup_{t\uparrow\infty} t^{-1} \log \mathbb{P}(K_t > \lfloor \alpha t \rfloor) = \limsup_{\alpha\uparrow\infty} \limsup_{t\uparrow\infty} t^{-1} \log \mathbb{P}\Bigg(\sum_{i=1}^{\lfloor \alpha t \rfloor} X_i < t\Bigg) = -\infty,$$

and hence, the error terms play no role on the exponential scale with rate $t$. Now, let $(t_n)_{n\geq 0}$ be a subsequence such that

$$\begin{aligned}&\limsup_{t\uparrow\infty} t^{-1} \log \sum_{m=0}^{\lfloor \alpha t \rfloor} \mathbb{P}\Bigg(\sum_{i=1}^{m} Y_i \geq t(a-\varepsilon), t(1-\varepsilon) \leq \sum_{i=1}^{m} X_i < t\Bigg)\\ &\quad= \lim_{n\uparrow\infty} t_n^{-1} \log \sum_{m=0}^{\lfloor \alpha t_n \rfloor} \mathbb{P}\Bigg(\sum_{i=1}^{m} Y_i \geq t_n(a-\varepsilon), t_n(1-\varepsilon) \leq \sum_{i=1}^{m} X_i < t_n\Bigg)\end{aligned}$$

and define

$$\beta_{n,\varepsilon} := \operatorname{argmax}\Bigg\{\mathbb{P}\Bigg(\sum_{i=1}^{\lfloor \beta t_n \rfloor} Y_i \geq t_n(a-\varepsilon), t_n(1-\varepsilon) \leq \sum_{i=1}^{\lfloor \beta t_n \rfloor} X_i < t_n\Bigg) : 0 \leq \beta \leq \alpha\Bigg\},$$

where we simply take $\beta_{n,\varepsilon}$ to be the smallest solution in case of ambiguity. Then, for a suitable further sub-sequence $(t_{n_k})_{k\geq 0}$, we have

$$\begin{aligned}&\lim_{n\uparrow\infty} t_n^{-1} \log \sum_{m=0}^{\lfloor \alpha t_n \rfloor} \mathbb{P}\Bigg(\sum_{i=1}^{m} Y_i \geq t_n(a-\varepsilon), t_n(1-\varepsilon) \leq \sum_{i=1}^{m} X_i < t_n\Bigg)\\ &\quad\leq \limsup_{n\uparrow\infty} t_n^{-1} \log \mathbb{P}\Bigg(\sum_{i=1}^{\lfloor \beta_{n,\varepsilon} t_n \rfloor} Y_i \geq t_n(a-\varepsilon), t_n(1-\varepsilon) \leq \sum_{i=1}^{\lfloor \beta_{n,\varepsilon} t_n \rfloor} X_i < t_n\Bigg)\\ &\quad= \lim_{k\uparrow\infty} t_{n_k}^{-1} \log \mathbb{P}\Bigg(\sum_{i=1}^{\lfloor \beta_{n_k,\varepsilon} t_{n_k} \rfloor} Y_i \geq t_{n_k}(a-\varepsilon), t_{n_k}(1-\varepsilon) \leq \sum_{i=1}^{\lfloor \beta_{n_k,\varepsilon} t_{n_k} \rfloor} X_i < t_{n_k}\Bigg),\end{aligned}$$

and we note that $(\beta_{n_k,\varepsilon})_{k\geq 0}$ must contain a convergent sub-sequence $(\beta_{n_{k_l},\varepsilon})_{l\geq 0}$ with limit $0 \leq \beta^*_\varepsilon \leq \alpha$. Hence,

$$\lim_{l\uparrow\infty} t_{n_{k_l}}^{-1} \log \mathbb{P}\left( \sum_{i=1}^{\lfloor \beta_{n_{k_l},\varepsilon} t_{n_{k_l}} \rfloor} Y_i \geq t_{n_{k_l}}(a-\varepsilon),\ t_{n_{k_l}}(1-\varepsilon) \leq \sum_{i=1}^{\lfloor \beta_{n_{k_l},\varepsilon} t_{n_{k_l}} \rfloor} X_i < t_{n_{k_l}} \right)$$
$$\leq -\inf\{\beta \mathcal{J}_{\theta,\lambda}\big(1/\beta, a/\beta\big) \colon \beta > 0\},$$

since $\limsup_{t\uparrow\infty} t^{-1} \log\lfloor \alpha t\rfloor = 0$ and we could use continuity.

***Step 3:*** *General sets.* Note that, by symmetry, $\mathbb{P}(\mathcal{Y}_t \geq at) = \mathbb{P}(\mathcal{Y}_t \leq -at)$ and, hence, for $O \subset \mathbb{R}^2$ open, we find that, for all $z = (x, y) \in O$, there exists $\varepsilon_z > 0$ such that $\bar{B}_{\varepsilon_z}(z) \subset O$. In particular, for all $z \in O$,

$$\mathbb{P}(t^{-1}\mathcal{Y}_t \in O) \geq \mathbb{P}(z - \varepsilon_z < t^{-1}\mathcal{Y}_t < z + \varepsilon_z) = \mathbb{P}(t^{-1}\mathcal{Y}_t > -z - \varepsilon_z) - \mathbb{P}(t^{-1}\mathcal{Y}_t \geq -z + \varepsilon_z),$$

where the rate to zero is dominated by the first summand. Taking an infimum over $z \in O$ gives the desired result. The upper bound can be proved similarly.

**Step 4:** *Scaling.* If we denote by $\mathcal{Y}_{\lambda,t}$ the vertical displacement at time $t$ in the navigation based on a Poisson point process with intensity $\lambda > 0$, then, by scaling both the coordinates by $\sqrt{\lambda}$, we find that $\mathcal{Y}_{\lambda,t}$ and $\mathcal{Y}_{1,\sqrt{\lambda}t}/\sqrt{\lambda}$ are equal in distribution. Therefore, for any $x > 0$,

$$\begin{aligned}
\mathcal{I}_{\lambda,\theta}(x) &= -\lim_{t\uparrow\infty} t^{-1} \log \mathbb{P}(t^{-1}\mathcal{Y}_{\lambda,t} > x) \\
&= -\lim_{t\uparrow\infty} \sqrt{\lambda}(\sqrt{\lambda}t)^{-1} \log \mathbb{P}((\sqrt{\lambda}t)^{-1}\mathcal{Y}_{1,\sqrt{\lambda}t} > x) \\
&= \sqrt{\lambda}\mathcal{I}_{1,\theta}(x).
\end{aligned}$$

By a similar calculation, we also find that, for any $x < 0$, $\mathcal{I}_{\lambda,\theta}(x) = \sqrt{\lambda}\mathcal{I}_{1,\theta}(x)$. This proves the theorem.

### 7.2. The dependent case $\pi/4 < \theta < \pi/2$

*Proof of Lemma* 3. The proof follows again from an application of the multivariate Cramér theorem for empirical means of sequences of i.i.d. random variables; see [7, Corollary 6.1.6]. We need to show existence of $\gamma = (\gamma_1, \gamma_2) \in \mathbb{R}^2 \setminus \{0\}$, such that $J'_{\lambda,\theta}(\gamma) < \infty$, but this is an immediate consequence of Lemma 2.

*Proof of Lemma* 4. Observe that if, for all $1 \leq m \leq n$, $\Phi_m \in (\pi/2 - \theta, \theta]$ then, for all $1 \leq m \leq n$, $H_m \neq \emptyset$, which then ensures that $\tau_1^\theta > n$. Therefore,

$$\mathbb{P}(\tau_1^\theta > n) \geq \mathbb{P}(\Phi_m \in (\pi/2 - \theta, \theta] \text{ for all } 1 \leq m \leq n).$$

Note that, conditioned on $V_{m-1}$, $R_m$ and $H_{m-1}$, $\Phi_m$ is uniformly distributed on the set

$$\{\varphi \in [-\theta, \theta] \colon V_{m-1} + (R_m, \varphi) \notin H_{m-1}\}.$$

Hence, with $\ell$ denoting the Lebesgue measure,

$$\mathbb{P}(\Phi_m \in (\pi/2 - \theta, \theta] | V_{m-1}, R_m, H_{m-1}) = \frac{\ell(\{\varphi \in (\pi/2 - \theta, \theta] \colon V_{m-1} + (R_m, \varphi) \notin H_{m-1}\})}{\ell(\{\varphi \in [-\theta, \theta] \colon V_{m-1} + (R_m, \varphi) \notin H_{m-1}\})}.$$

Also, note that if, for all $1 \le i < m$, $\Phi_i > 0$ then $\mathcal{C}^u_\theta(V_{m-1}) \cap H_{m-1} = \emptyset$. Therefore, in this case,

$$\frac{\ell(\{\varphi \in (\pi/2-\theta, \theta] \colon V_{m-1} + (R_m, \varphi) \notin H_{m-1}\})}{\ell(\{\varphi \in [-\theta, \theta] \colon V_{m-1} + (R_m, \varphi) \notin H_{m-1}\})} \ge \frac{\ell((\pi/2-\theta, \theta])}{\ell([-\theta, \theta])} = \frac{2\theta - \pi/2}{2\theta} = \frac{4\theta - \pi}{4\theta}.$$

As a result, we find that

$$\mathbb{P}(\Phi_m \in (\pi/2-\theta, \theta] | \Phi_i \in (\pi/2-\theta, \theta] \text{ for all } 1 \le i < m) \ge (4\theta - \pi)/(4\theta),$$

and thus,

$$\begin{aligned}\mathbb{P}(\tau > n) &\ge \mathbb{P}(\Phi_m \in (\pi/2-\theta, \theta] \text{ for all } 1 \le m \le n)\\ &= \prod_{m=1}^{n} \mathbb{P}(\Phi_m \in (\pi/2-\theta, \theta] | \Phi_i \in (\pi/2-\theta, \theta] \text{ for all } 1 \le i < m)\\ &\ge ((4\theta - \pi)/(4\theta))^n,\end{aligned}$$

which completes the proof.

*Proof of Theorem* 3. We proceed via several steps. First, note that, for all $\delta > 0$,

$$\limsup_{t \uparrow \infty} t^{-1} \log \mathbb{P}\big(|\mathcal{Y}_t - \hat{\mathcal{Y}}'_t| \ge \delta t\big) = -\infty,$$

where $\hat{\mathcal{Y}}'_t := \sum_{i=1}^{K_t} Y_i$, since

$$\begin{aligned}\mathbb{P}\big(|\mathcal{Y}_t - \hat{\mathcal{Y}}'_t| \ge \delta t\big) &\le \mathbb{P}(|Y_{K_t+1}| > \delta t)\\ &\le \mathbb{P}(K_t \ge \lfloor t \rfloor^2) + \mathbb{P}\left(\max_{i=1}^{\lfloor t \rfloor^2} |Y_i| > \delta t\right)\\ &\le \mathbb{P}\left(\sum_{i=1}^{\lfloor t \rfloor^2} X_i < t\right) + \mathbb{P}\left(\max_{i=1}^{\lfloor t \rfloor^2} |Y_i| > \delta t\right)\\ &\le \mathbb{P}\left(\sum_{i=1}^{\lfloor t \rfloor^2} \underline{R}_i \cos\theta < t\right) + \mathbb{P}\left(\max_{i=1}^{\lfloor t \rfloor^2} \overline{R}_i > \delta t\right)\\ &\le \mathrm{e}^t \mathbb{E}[\mathrm{e}^{-\underline{R}\cos\theta}]^{\lfloor t \rfloor^2} + t^2 \mathrm{e}^{-\lambda(\theta \wedge (\pi/2-\theta))\delta^2 t^2},\end{aligned}$$

and hence, for the large deviations, we can focus on $\hat{\mathcal{Y}}'_t$. For this, as in (3), we can further split $\hat{\mathcal{Y}}'_t$ into

$$\hat{\mathcal{Y}}'_t = \sum_{i=1}^{K'_t} Y'_i + \sum_{j=\tau^\theta_{K'_t}+1}^{K_t} Y_j =: \mathcal{Y}'_t + \hat{Y}_t.$$

The first summand is a sum of i.i.d. segments, but the second summand is a sum of dependent random variables. The challenge comes from the fact that all the involved random variables $Y'_i$, $\hat{Y}_t$ have exponential tails and, therefore, contribute on the large-deviation scale $t$.

***Step 1:*** *Upper bound for upper tail.* For all $\alpha > 0$, we can bound

$$\mathbb{P}(\hat{\mathcal{Y}}'_t \geq at) \leq \mathbb{P}\left(\sum_{i=1}^{K'_t} Y'_i + \hat{Y}_t \geq at,\, K'_t \leq \lfloor \alpha t \rfloor\right) + \mathbb{P}(K'_t > \lfloor \alpha t \rfloor)$$

$$\leq \sup\nolimits_{0<\beta\leq\alpha} \alpha t \mathbb{P}\left(\sum_{i=1}^{\lfloor \beta t \rfloor} Y'_i + \hat{Y}_t \geq at,\, K'_t = \lfloor \beta t \rfloor\right) + \mathbb{P}(K'_t > \lfloor \alpha t \rfloor),$$

where the second summand plays no role on the exponential scale with rate $t$ since

$$\limsup_{\alpha\uparrow\infty} \limsup_{t\uparrow\infty} t^{-1} \log \mathbb{P}(K'_t > \lfloor \alpha t \rfloor) = \limsup_{\alpha\uparrow\infty} \limsup_{t\uparrow\infty} t^{-1} \log \mathbb{P}\left(\sum_{i=1}^{\lfloor \alpha t \rfloor} X'_i < t\right) = -\infty.$$

Similarly, for all $0 < \beta \leq \alpha$,

$$\limsup_{\gamma\uparrow\infty} \limsup_{t\uparrow\infty} t^{-1} \log \mathbb{P}\left(\sum_{i=1}^{\lfloor \beta t \rfloor} Y'_i > \gamma t\right) = -\infty$$

and

$$\limsup_{\delta\uparrow\infty} \limsup_{t\uparrow\infty} t^{-1} \log \mathbb{P}(\tau_1^\theta > \delta t) = -\infty,$$

and hence, it suffices to further bound as

$$\mathbb{P}\left(\sum_{i=1}^{\lfloor \beta t \rfloor} Y'_i + \hat{Y}_t \geq at,\, K'_t = \lfloor \beta t \rfloor\right)$$

$$\leq \sup\nolimits_{|b|<\gamma,\ c<1} \gamma t^2 \mathbb{P}\left(\sum_{i=1}^{\lfloor \beta t \rfloor} Y'_i \mathrel{\hat{=}} bt,\, \sum_{i=1}^{\lfloor \beta t \rfloor} X'_i \mathrel{\hat{=}} ct\right)$$

$$\times \sup\nolimits_{d\leq\delta} \delta t \mathbb{P}\left(\sum_{i=1}^{\lfloor dt \rfloor} Y_i \geq (a-b)t,\, \sum_{i=1}^{\lfloor dt \rfloor} X_i < (1-c)t,\, \sum_{i=1}^{\lfloor dt \rfloor+1} X_i \geq (1-c)t,\, \tau_1^\theta > dt\right)$$

$$+ \mathbb{P}\left(\sum_{i=1}^{\lfloor \beta t \rfloor} Y'_i > \gamma t\right) + \mathbb{P}(\tau_1^\theta > \delta t),$$

where we write $x \mathrel{\hat{=}} y$ if and only if $\lfloor x \rfloor = \lfloor y \rfloor$ and use independence. Now, we can combine Assumption 1 and Lemma 3 to obtain

$$\limsup_{t\uparrow\infty} t^{-1} \log \mathbb{P}(\hat{\mathcal{Y}}'_t \geq at)$$

$$\leq -\inf_{b\in\mathbb{R},\ c\in(0,1)} \Big\{ \inf\{\beta \mathcal{J}'_{\lambda,\theta}(c/\beta, b/\beta)\colon \beta > 0\} + \inf\{d\mathcal{H}_{\lambda,\theta}((1-c)/d, (a-b)/d)\colon d > 0\}\Big\},$$

as desired.

***Step 2:*** *Lower bound for upper tail.* Next, we consider the lower bound for the upper tail with $a > 0$. Then, for $b \in \mathbb{R}$, $0 < c < 1$ and $\beta, d > 0$,

$$\begin{aligned}\mathbb{P}(\hat{\mathcal{Y}}'_t > at) &\geq \mathbb{P}\left(\sum_{i=1}^{K'_t} Y'_i > bt,\, K'_t = \lfloor \beta t \rfloor,\, \hat{Y}_t > (a-b)t\right)\\ &\geq \mathbb{P}\left(\sum_{i=1}^{\lfloor \beta t \rfloor} Y'_i > bt,\, \sum_{i=1}^{\lfloor \beta t \rfloor} X'_i \hat{=} ct\right)\\ &\quad \times \mathbb{P}\left(\sum_{i=1}^{\lfloor dt \rfloor} Y_i > (a-b)t,\, \sum_{i=1}^{\lfloor dt \rfloor} X_i < (1-c)t,\, \sum_{i=1}^{\lfloor dt \rfloor + 1} X_i \geq (1-c)t,\, \tau_1^\theta > \lfloor dt \rfloor\right),\end{aligned}$$

where we used independence. Consequently,

$$\liminf_{t\uparrow\infty} t^{-1} \log \mathbb{P}(\hat{\mathcal{Y}}'_t > at) \geq -\beta \mathcal{I}'_{\lambda,\theta}(c/\beta, b/\beta) - d\mathcal{H}_{\lambda,\theta}((1-c)/d, (a-b)/d).$$

Optimizing first with respect to $\beta$ and $d$ in the individual summands and then with respect to $b$, $c$ in the joint expression, we arrive at the desired lower bound that matches the upper bound.

***Step 3:*** *General sets.* Using symmetry and the previous steps, we can follow the exact same arguments as in the independent case, step 4 in the proof of Theorem 2, to arrive at the large-deviation principle.

***Step 4:*** *Scaling.* A scaling argument similar to that in the proof of Theorem 2 implies that $\mathcal{I}'_{\lambda,\theta}(x) = \sqrt{\lambda}\mathcal{I}'_{1,\theta}(x)$ for all $x \in \mathbb{R}$.

## Appendix A. Details for Remark 1

Let us finally present some more details for Remark 1. Let $0 < \theta \leq \pi/4$ and note that $\gamma \mapsto J_{\lambda,\theta}$ is strictly convex and twice differentiable. Hence, fixing $x$, with $|x| < \tan\theta$, and writing $\zeta = 1/\beta$, for every $\zeta > 0$, there exists a unique $\gamma(\zeta)$ such that

$$\mathcal{J}(\zeta) := \mathcal{J}_{\lambda,\theta}(\zeta, \zeta x) = \langle \gamma(\zeta), (\zeta, \zeta x)\rangle - J_{\lambda,\theta}(\gamma(\zeta)),$$

where $\gamma(\zeta)$ is twice differentiable, due to the explicit-function theorem, and satisfies

$$\zeta = \frac{\mathbb{E}[\exp(\langle \gamma(\zeta), U\rangle)X]}{\mathbb{E}[\exp(\langle \gamma(\zeta), U\rangle)]} \quad \text{and} \quad \zeta x = \frac{\mathbb{E}[\exp(\langle \gamma(\zeta), U\rangle)Y]}{\mathbb{E}[\exp(\langle \gamma(\zeta), U\rangle)]}, \tag{A.1}$$

where we suppressed the dependence on $\lambda$ and $\theta$ in the notation and let $U = (X, Y)$ represent the random variable for the first step. In particular, the minimizing of $\zeta = \zeta(x)$ in the definition of $\mathcal{I}_{\lambda,\theta}(x)$ satisfies

$$0 = -\frac{1}{\zeta^2}\mathcal{J}(\zeta) + \frac{1}{\zeta}\dot{\mathcal{J}}(\zeta) = -\frac{1}{\zeta^2}\mathcal{J}(\zeta) + \frac{1}{\zeta}(\gamma_1(\zeta) + x\gamma_2(\zeta)),$$

where we used (A.1). In other words, $\mathcal{J}(\zeta) = \zeta\gamma_1(\zeta) + \zeta x\gamma_2(\zeta) = \langle \gamma(\zeta), (\zeta, \zeta x)\rangle$, which implies that $\mathbb{E}[\exp(\langle \gamma(\zeta(x)), U\rangle)] = 1$ and, moreover,

$$\mathcal{I}_{\lambda,\theta}(x) = \varrho_1(x) + x\varrho_2(x),$$

where we abbreviated $\varrho_1(x) := \gamma_1(\zeta(x))$ and $\varrho_2(x) := \gamma_2(\zeta(x))$. This in particular provides us with a representation of the large-deviation rate function in terms of the unique solution to (A.1). Hence,

$$\ddot{\mathcal{I}}_{\lambda,\theta}(0) = \ddot{\varrho_1}(0) + 2\dot{\varrho_2}(0).$$

Now, in order to analyze this, note that, by symmetry, $\mathcal{I}_{\lambda,\theta}(0) = 0$ implies that $\varrho_1(0) = \varrho_2(0) = 0$ and, hence, using the left-hand side of (A.1), $\zeta(0) = \mathbb{E}[X]$. Also, note that $\dot{\mathcal{I}}_{\lambda,\theta}(0) = 0$ implies that $\dot{\varrho}_1(0) = 0$. Considering the first derivative on the left-hand side of (A.1), we have

$$\dot{\zeta}(x) = \frac{\mathbb{E}[\exp(\langle\varrho(x), U\rangle)X\langle\dot{\varrho}(x), U\rangle]}{\mathbb{E}[\exp(\langle\varrho(x), U\rangle)]} - \frac{\mathbb{E}[\exp(\langle\varrho(x), U\rangle)X]\mathbb{E}[\exp(\langle\varrho(x), U\rangle)\langle\dot{\varrho}(x), U\rangle]}{\mathbb{E}[\exp(\langle\varrho(x), U\rangle)]^2},$$

which implies that $\dot{\zeta}(0) = \mathbb{E}[X\langle\dot{\varrho}(0), U\rangle] - \mathbb{E}[X]\mathbb{E}[\langle\dot{\varrho}(0), U\rangle] = \dot{\varrho}_1(0)\mathbb{V}[X] = 0$. By symmetry, it is immediate that, for the right-hand side of (A.1),

$$\begin{aligned}&\dot{\zeta}(x)x + \zeta(x)\\ &\quad= \frac{\mathbb{E}[\exp(\langle\varrho(x), U\rangle)Y\langle\dot{\varrho}(x), U\rangle]}{\mathbb{E}[\exp(\langle\varrho(x), U\rangle)]} - \frac{\mathbb{E}[\exp(\langle\varrho(x), U\rangle)Y]\mathbb{E}[\exp(\langle\varrho(x), U\rangle)\langle\dot{\varrho}(x), U\rangle]}{\mathbb{E}[\exp(\langle\varrho(x), U\rangle)]^2},\end{aligned}$$

and hence, $\mathbb{E}[X] = \zeta(0) = \mathbb{E}[Y\langle\dot{\varrho}(0), U\rangle] - \mathbb{E}[Y]\mathbb{E}[\langle\dot{\varrho}(0), U\rangle] = \dot{\varrho_2}(0)\mathbb{V}[Y]$, that is, $\dot{\varrho_2}(0) = \mathbb{E}[X]/\mathbb{E}[Y^2]$. Furthermore, considering the left-hand side of (A.1), we have

$$\begin{aligned}\ddot{\zeta}(x) &= \frac{\mathbb{E}[\exp(\langle\varrho(x), U\rangle)X\langle\dot{\varrho}(x), U\rangle^2] + \mathbb{E}[\exp(\langle\varrho(x), U\rangle)X\langle\ddot{\varrho}(x), U\rangle]}{\mathbb{E}[\exp(\langle\varrho(x), U\rangle)]}\\ &\quad - 2\frac{\mathbb{E}[\exp(\langle\varrho(x), U\rangle)X\langle\dot{\varrho}(x), U\rangle]\mathbb{E}[\exp(\langle\varrho(x), U\rangle)\langle\dot{\varrho}(x), U\rangle]}{\mathbb{E}[\exp(\langle\varrho(x), U\rangle)]^2}\\ &\quad - \frac{\mathbb{E}[\exp(\langle\varrho(x), U\rangle)X]\big(\mathbb{E}[\exp(\langle\varrho(x), U\rangle)\langle\dot{\varrho}(x), U\rangle^2] + \mathbb{E}[\exp(\langle\varrho(x), U\rangle)\langle\ddot{\varrho}(x), U\rangle]\big)}{\mathbb{E}[\exp(\langle\varrho(x), U\rangle)]^2}\\ &\quad + 2\frac{\mathbb{E}[\exp(\langle\varrho(x), U\rangle)X]\mathbb{E}[\exp(\langle\varrho(x), U\rangle)\langle\dot{\varrho}(x), U\rangle]^2\mathbb{E}[\exp(\langle\varrho(x), U\rangle)]}{\mathbb{E}[\exp(\langle\varrho(x), U\rangle)]^4},\end{aligned}$$

which gives a representation of $\ddot{\varrho_1}(0)$ in terms of mixed moments,

$$\begin{aligned}\ddot{\zeta}(0) &= \mathbb{E}[X\langle\dot{\varrho}(0), U\rangle^2] + \mathbb{E}[X\langle\ddot{\varrho}(0), U\rangle] - 2\mathbb{E}[X\langle\dot{\varrho}(0), U\rangle]\mathbb{E}[\langle\dot{\varrho}(0), U\rangle]\\ &\quad - \mathbb{E}[X]\big(\mathbb{E}[\langle\dot{\varrho}(0), U\rangle^2] + \mathbb{E}[\langle\ddot{\varrho}(0), U\rangle]\big) + 2\mathbb{E}[X]\mathbb{E}[\langle\dot{\varrho}(0), U\rangle]^2\\ &= \dot{\varrho_2}(0)^2\mathbb{E}[XY^2] + \ddot{\varrho_1}(0)\mathbb{E}[X^2] - \mathbb{E}[X]\big(\dot{\varrho_2}(0)^2\mathbb{E}[Y^2] + \ddot{\varrho_1}(0)\mathbb{E}[X]\big)\\ &= \dot{\varrho_2}(0)^2(\mathbb{E}[XY^2] - \mathbb{E}[X]\mathbb{E}[Y^2]) + \ddot{\varrho_1}(0)\mathbb{V}[X].\end{aligned}$$

Putting things together, we get the representation

$$\ddot{\mathcal{I}}_{\lambda,\theta}(0) = \left(\ddot{\zeta}(0) - \frac{(\mathbb{E}[XY^2] - \mathbb{E}[X]\mathbb{E}[Y^2])\mathbb{E}[X]^2}{\mathbb{E}[Y^2]^2}\right)/\mathbb{V}[X] + 2\mathbb{E}[X]/\mathbb{E}[Y^2].$$

In order to see that indeed $\rho_{\lambda,\theta} \neq \ddot{\mathcal{I}}_{\lambda,\theta}(0)/2$, first note that $\ddot{\zeta}(0) \leq 0$ since, for a deviation event in the vertical direction, it is not beneficial to increase the number of horizontal steps. This would increase the averaging effect and would make it harder to achieve the vertical deviation. Hence,

$$\frac{\ddot{\mathcal{I}}_{\lambda,\theta}(0)}{2} \leq \frac{\mathbb{E}[X]}{\mathbb{E}[Y^2]} - \frac{(\mathbb{E}[XY^2] - \mathbb{E}[X]\mathbb{E}[Y^2])\mathbb{E}[X]^2}{2\mathbb{E}[Y^2]^2\mathbb{V}[X]}$$
$$= \rho_{\lambda,\theta}\left(2 - \frac{(\mathbb{E}[XY^2] - \mathbb{E}[X]\mathbb{E}[Y^2])\mathbb{E}[X]}{\mathbb{E}[Y^2]\mathbb{V}[X]}\right),$$

and it suffices to show that

$$\mathrm{Cov}(X, Y^2)\mathbb{E}[X] > \mathbb{E}[Y^2]\mathbb{V}[X].$$

However, using the moment-generating function (24), we can compute the left- and right-hand side of this, for example, for $\theta = \pi/4$ and $\lambda = 2$, numerically to see that it can be satisfied.

## Acknowledgements

The authors would like to thank Christian Hirsch, Kumarjit Saha, Wolfgang König, and Alexander Zass for many fruitful discussions about the topic and for mentioning important references. The author also thanks the anonymous referees for their careful reading and detailed comments, which have helped to improve the paper.

## Funding information

This work is supported by the Leibniz Association within the Leibniz Junior Research Group on *Probabilistic Methods for Dynamic Communication Networks* as part of the Leibniz Competition (grant no. J105/2020). The third author would also like to thank the grant *ERC NEMO* (grant no. 788851) of the research group DYOGENE at INRIA Paris.

## Competing interests

There were no competing interests to declare which arose during the preparation or publication process of this paper.

## References

[1] ARCONES, M. A. (2003). Large deviations of empirical processes. In *High Dimensional Probability, III (Sandjberg, 2002)*, Vol. 55 of Progress in Probability. Birkhäuser, Basel, pp. 205–223.

[2] ASMUSSEN, S. (2003). *Applied Probability and Queues*, 2nd edn., Vol. 51 of Applications of Mathematics (New York), Stochastic Modelling and Applied Probability. Springer-Verlag, New York.

[3] BACCELLI, F. AND BORDENAVE, C. (2007). The radial spanning tree of a Poisson point process. *Ann. Appl. Probab.* **17**(1), 305–359.

[4] BONICHON, N. AND MARCKERT, J. F. (2011). Asymptotics of geometrical navigation on a random set of points in the plane. *Adv. Appl. Probab.* **43**(4), 899–942.

[5] BORDENAVE, C. (2008). Navigation on a Poisson point process. *Ann. Appl. Probab.* **18**(2), 708–746.

[6] COUPIER, D., SAHA, K., SARKAR, A. AND TRAN, V. C. (2021). The 2d-directed spanning forest converges to the Brownian web. *Ann. Probab.* **49**(1), 435–484.

[7] DEMBO, A. AND ZEITOUNI, O. (2010). *Large Deviations Techniques and Applications*, Vol. 38 of Stochastic Modelling and Applied Probability. Springer-Verlag, Berlin. Corrected reprint of the second (1998) edition.

[8] HOLLANDER, F. den. (2000). *Large Deviations*, Vol. 14 of Fields Institute Monographs. American Mathematical Society, Providence, RI.

[9] EICHELSBACHER, P. AND LÖWE, M. (2003). Moderate deviations for i.i.d. random variables. *ESAIM Probab. Stat.* **7**, 209–218.

[10] HIRSCH, C., JAHNEL, J., KEELER, P. AND PATTERSON, R. I. A. (2017). Traffic flow densities in large transport networks. *Adv. Appl. Probab.* **49**(4), 1091–1115.

[11] HOWARD, C. D. AND NEWMAN, C. M. (1997). Euclidean models of first-passage percolation. *Probab. Theory Related Fields* **108**(2), 153–170.

[12] JAHNEL, B. AND KÖNIG, W. (2020). *Probabilistic Methods in Telecommunications*. *Compact Textbooks in Mathematics*. Birkhäuser, Cham.

[13] LEDOUX, M. (1992). Sur les déviations modérées des sommes de variables aléatoires vectorielles indépendantes de même loi. *Ann. Inst. H. Poincaré Probab. Statist.* **28**(2), 267–280.

[14] ROCKAFELLAR, R. T. (1967). Conjugates and Legendre transforms of convex functions. *Canad. J. Math.* **19**, 200–205.

[15] ROY, R., SAHA, K. AND SARKAR, A. (2016) Random directed forest and the Brownian web. *Ann. Inst. Henri Poincaré Probab. Stat.* **52**(3), 1106–1143.

[16] ROY, R., SAHA, K. AND SARKAR, A. (2023). Scaling limit of a drainage network model on perturbed lattice. arXiv:2302.09489.

[17] VARADHAN, S. R. S. (1984). *Large Deviations and Applications*, Vol. 46 of CBMS-NSF Regional Conference Series in Applied Mathematics. Society for Industrial and Applied Mathematics (SIAM), Philadelphia, PA.

[18] BALL, K. AND CHAIN, H. (1988). *Kurtosis: A Critical Review*, 2nd edn. John Wiley, New York.

[19] BOYD, W. (1978). Hyperbolic distributions. Doctoral Thesis, University of Boston School of Mathematics.

[20] SICHEL, H. S., KLEINGELD, W. J. AND ASSIBEY-BONSU, W. (1992). A comparative study of three frequency-distribution models for use in ore valuation. *J. S. Afr. Inst. Min. Met.* **92**, 91–99.